\documentclass[12pt,letter]{article}
\usepackage{graphicx}
\usepackage{epstopdf}
\usepackage{subfig}
\usepackage{fullpage}
\usepackage{amsmath, amssymb}
\usepackage{color, colortbl}
\usepackage{xcolor}
\usepackage{float}
\usepackage{url}

\title{Multigrid-in-time for sensitivity analysis of chaotic dynamical systems}
\author{Patrick Blonigan, MIT and Qiqi Wang, MIT}

\begin{document}

\maketitle

\section*{Abstract}

The following paper discusses the application of a multigrid-in-time scheme to Least Squares Shadowing (LSS), a novel sensitivity analysis method for chaotic dynamical systems.  While traditional sensitivity analysis methods break down for chaotic dynamical systems, LSS is able to compute accurate gradients.  Multigrid is used because LSS requires solving a very large Karush-–Kuhn-–Tucker (KKT) system constructed from the solution of the dynamical system over the entire time interval of interest.

Several different multigrid-in-time schemes are examined, and a number of factors were found to heavily influence the convergence rate of multigrid-in-time for LSS.  These include the iterative method used for the smoother, how the coarse grid system is formed and how the least squares objective function at the center of LSS is weighted.  

\section{Introduction}
\label{c:intro}

Although computers are becoming more powerful at a rapid pace, the size of the problems we wish to solve are growing even faster.  An example of a potentially very large problem is Least Squares Shadowing (LSS), a novel computational method that can be used for sensitivity analysis of chaotic dynamical systems. 

LSS requires solving a linear system of $\mathcal{O}(mn)$ equations, where $m$ is the number of time steps and $n$ is the number of dimensions or degrees of freedom of the system, both of which can be very large ($\sim 10^5$ or greater) for many problems of interest in a number of fields including fluid dynamics.  Fortunately this system can be written as a 2nd order boundary value problem in time with homogeneous Dirichlet boundary conditions. Therefore a multigrid-in-time scheme is attractive because of its fast convergence relative to other iterative methods for many boundary value problems in time \cite{Briggs:2000:MGtutorial}. 

\subsection{Motivation for Least Squares Shadowing}

Sensitivity analysis of systems governed by ordinary differential
equations (ODEs) and partial differential equations (PDEs) is important in many
fields of science and engineering.  Its goal is to compute sensitivity
derivatives of key quantities of interest to parameters that influence
the system.  Applications of sensitivity analysis in science and
engineering include design optimization, inverse problems,
data assimilation, and uncertainty quantification.

Adjoint based sensitivity analysis is especially powerful in many
applications, due to its efficiency when the number of parameters is
large.  In aircraft design, for example, the number of geometric parameters that
define the aerodynamic shape is very large.  As a result, the adjoint
method of sensitivity analysis has proven to be very successful for aircraft design
 \cite{Jameson:1988:adj}, \cite{Reuther:2001:adjAC}.  
Similarly, the
adjoint method has been an essential tool for adaptive grid methods for solving PDEs \cite{Darmofal:2002:adapt}, solving inverse problems in
seismology, and for assimilating observation data for weather
forecasting.  

Sensitivity analysis for chaotic dynamical systems is important because of the prevalence of chaos in many scientific and engineering fields.  One example is highly turbulent gas flow of mixing and combustion processes in jet engines.  In this example, and in other applications with periodic or chaotic characteristics, statistical averaged quantities such as mean temperature and mean aerodynamic forces are of interest.  Therefore, the general problem we seek to solve with sensitivity analysis is:
 
\begin{equation}
\mbox{Given } \frac{du}{dt} = f(u, \xi), \quad
\overline{J} = \lim_{T\rightarrow\infty} \frac1T \int_0^T J(u,\xi)
dt,\quad
\mbox{Compute } \frac{\partial \overline{J}}{\partial \xi}
\label{e:gen_problem}
\end{equation}

\noindent Where $u$ is the state vector, $\xi$ is some parameter in the governing equation, $ \frac{du}{dt} = f(u, \xi)$ and $J(u,\xi)$ is some quantity of interest.  

Sensitivity analysis for chaotic dynamical systems is difficult because of the high sensitivity of these systems to the initial condition, known as the "Butterfly Effect".  Slightly different initial conditions will result in very different solutions, which diverge exponentially with time \cite{Lorenz:1963:det}.  This also results in exponential growth of sensitivities and therefore the sensitivity of long-time averaged quantities is not equal to the long-time average sensitivities of chaotic systems \cite{Lea:2000:climate_sens}.  Because the derivative and long-time average do not commute, the traditional adjoint method computes sensitivities that grow exponentially with time, as shown in the work done by Lea et al. \cite{Lea:2000:climate_sens}.  

The least squares shadowing (LSS) method does not encounter the exponential growth of sensitivities observed in traditional methods \cite{Wang:2013:LSS1}. The method requires that the chaotic system is {\em ergodic}; that long time behavior of the system is independent of initial conditions.  LSS finds a perturbed trajectory that does not diverge exponentially from some trajectory in phase space.  This non-diverging trajectory, called a ``shadow trajectory'', has its existence guaranteed by the shadowing lemma \cite{Pilyugin:1999:shadow} for a large number of chaotic systems and can be used to compute sensitivities.  The shadow trajectory is found by solving a quadratic programing (QP) problem with linear constraints.  It is the large size and bandwidth of the KKT system associated with the QP system that requires an efficient, iterative linear solver.  

\subsection{Outline}

The following paper discusses the implementation of LSS using a multigrid-in-time scheme to solve the KKT system derived from the QP problem. Several schemes were considered, each resulting in different convergence rates and computational costs.  The paper is structured as follows: first, LSS is presented through a discussion of non-linear dynamics and demonstrated on the Lorenz system.  Next, several multigrid-in-time schemes are presented and their performance solving for sensitivities of the Lorenz system is presented and discussed.  Finally, some concluding remarks and directions for future work are presented.  

\section{Least Squares Shadowing Method}
\label{c:LSS}

\subsection{Lyapunov Exponents and the Shadowing Lemma}

Before LSS and the shadowing lemma from which it is derived is discussed, a more in-depth discussion of Lyapunov exponents and Lyapunov covariant vectors is required.  For some system $\frac{du}{dt} = f(u)$, there exist Lyapunov covariant vectors $\phi_1(u),\phi_2(u),...,\phi_i(u)$ corresponding to each Lyapunov exponent $\Lambda_i$, which satisfy the equation \cite{Ginelli:2007:Lya}:
\[
\frac{d}{dt}\phi_i(u(t)) = \frac{\partial f}{\partial u} \cdot \phi_i(u(t)) - \Lambda_i \phi_i(u(t))
\]

\begin{figure}
\centering
\includegraphics[scale=0.7]{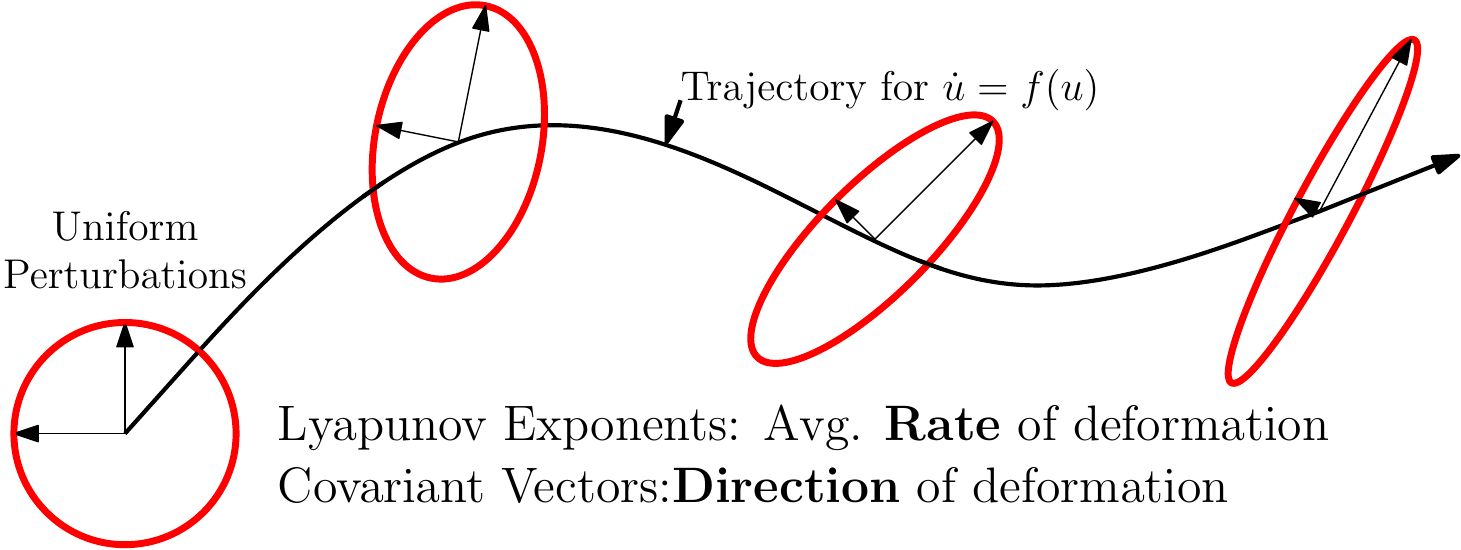}
\caption{Schematic of Lyapunov exponents and covariant vectors}
\label{f:Lyapunov}
\end{figure}

To understand what $\Lambda_i$ and $\phi_i$ represent, consider a sphere comprised of perturbations $\delta f$ to a system $\frac{du}{dt} = f(u)$ at some time, as shown in the far left of figure \ref{f:Lyapunov}. As this system evolves in time, this sphere expands in some directions, contracts in some, and remains unchanged in others.  The rate at which the sphere expands or contracts corresponds to the Lyapunov exponent $\Lambda_i$ and the corresponding direction of expansion or contraction is the Lyapunov covariant vector $\phi_i$.  It is important to note that the vectors $\phi_i$ are not the same as local Jacobian eigenvectors along a trajectory.  The $\phi_i$ vectors depend on all Jacobian eigenvectors along a trajectory. Also, the $\phi_i$ are not necessarily orthogonal, but the number of Lyapunov covariant vectors is the same as the number of dimensions of the system.  

\begin{figure}
\centering
\includegraphics[scale=0.3]{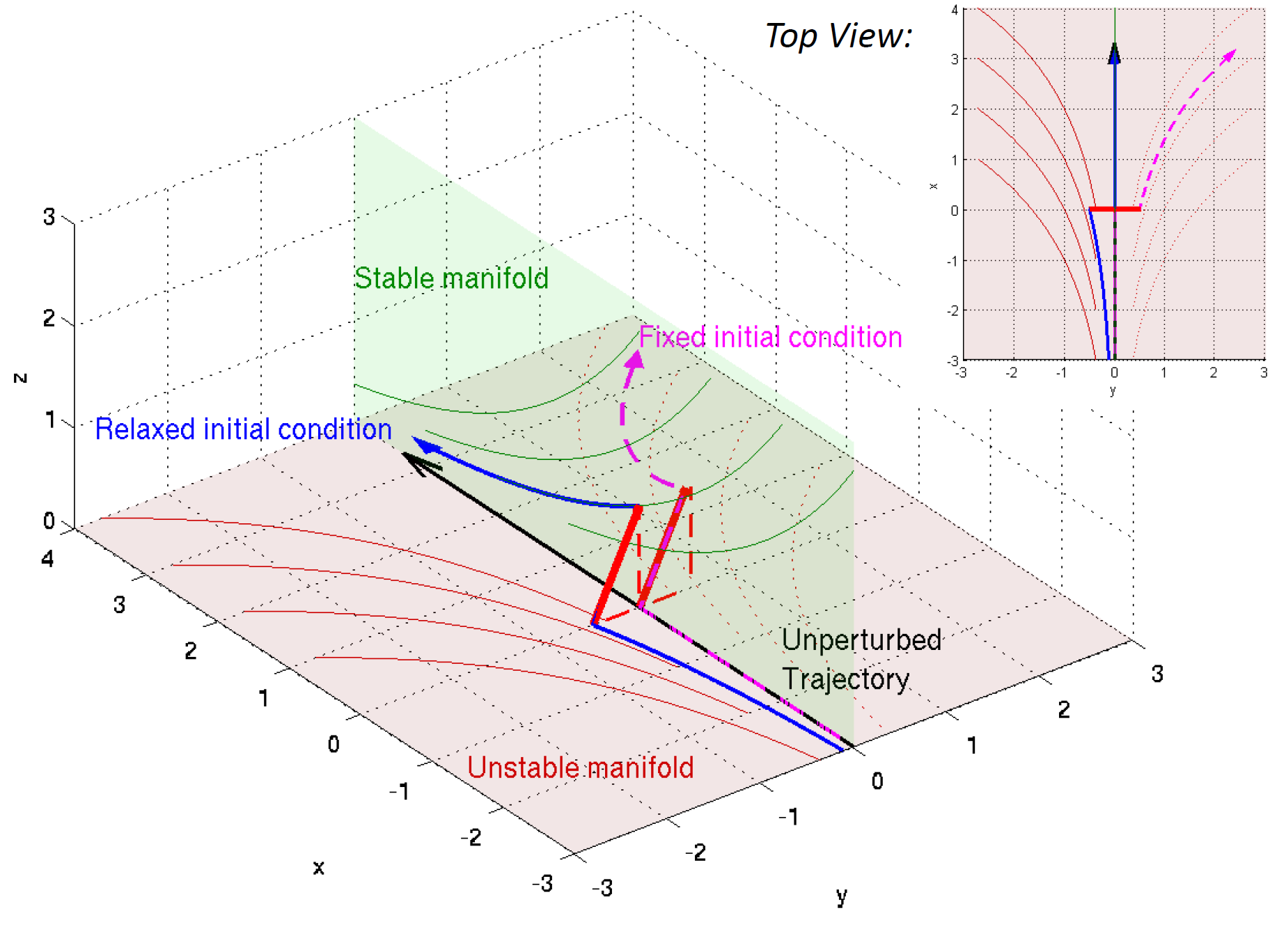}
\caption{Phase space trajectory of a chaotic dynamical system.  The unstable manifold, in red, is the space of all Lyapunov covariant vectors corresponding to positive exponents. The stable manifold, in green, corresponds to the space of all covariant vectors associated with negative exponents.  A perturbation to the system (in red) has components in both manifolds, and the unstable component causes the perturbed trajectory (pink) to diverge exponentially from the unperturbed trajectory (in black).  LSS chooses a perturbed trajectory with a different initial condition (in blue) that does not diverge from the unperturbed trajectory.  }
\label{f:LSS_schematic}
\end{figure}

A strange attractor, the type of attractor associated with chaotic dynamical systems, has at least one positive and one zero Lyapunov exponent \cite{Ginelli:2007:Lya}.  To illustrate the effect of the positive exponent, we consider figure \ref{f:LSS_schematic}.  We see that if the perturbed trajectory has the same initial condition as the unperturbed trajectory, the two trajectories diverge exponentially, leading to the issues with traditional sensitivity analysis discussed in the introduction.

However, the assumption of ergodicity means that it is not necessary to compare a perturbed and an unperturbed trajectory with the same initial condition.  Therefore, an initial condition can be chosen such that the perturbed and unperturbed trajectories do not diverge, resulting in the blue trajectory in figure \ref{f:LSS_schematic}. The existence of this trajectory, called a ``shadow trajectory'', follows from the shadowing lemma \cite{Pilyugin:1999:shadow} :
{\it
Consider a reference solution $u_r$ to

\[
\frac{du}{dt} = f(u,\xi)
\]

If this system has a Hyperbolic strange attractor and if some system parameter $\xi$ is slightly perturbed: 

For any $\delta > 0$ there exists $\varepsilon > 0$, such that for every $u_r$ that satisfies $\|d u_r / d t - f(u_r)\|<\varepsilon$, $0 \le t \le T$, there exists a true solution $u_s$ and a time transformation $\tau(t)$, such that $\|u_s(\tau(t))-u_r(t)\| < \delta$, $|1-d\tau/dt|<\delta$ and $d u_s/d \tau - f(u_s) = 0$, $0 \le \tau \le \mathcal{T}$.  

Note that $\|\cdot\|$ refers to distance in phase space
}

Therefore, relaxing the initial condition allows us to find the shadow trajectory $u_s(\tau)$.  The key assumption of the shadowing lemma is that the attractor associated with the system of interest is {\em hyperbolic}.  The key property of hyperbolic attractors for the shadowing lemma is that tangent space can be decomposed into stable, neutrally stable and unstable components everywhere on the attractor \cite{Hasselblatt:2002:hyperbolic}.  Another way to state this property is that the Lyapunov covariant vectors make up a basis for phase space at all points on the attractor.  Although this not the case for many attractors, including the Lorenz attractor, there is a {\em Chaotic Hypothesis} which states that many high-dimensional chaotic systems will behave as if they were hyperbolic \cite{Eyink:2004:ensmbl}.  For example, since the single point on the Lorenz attractor that is not hyperbolic is the unstable fixed point at the origin, most phase space trajectories do not pass through it and the shadowing lemma holds.  

\begin{figure}
\centering
\includegraphics[scale=0.6]{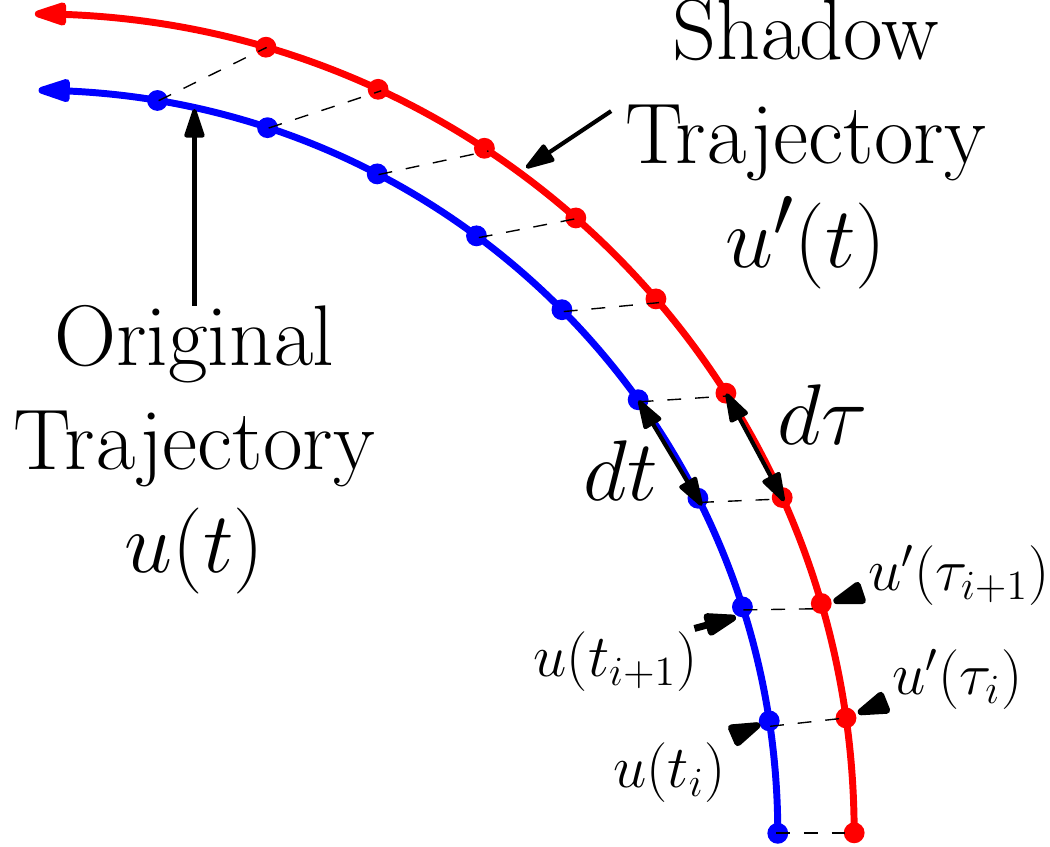}
\includegraphics[scale=0.6]{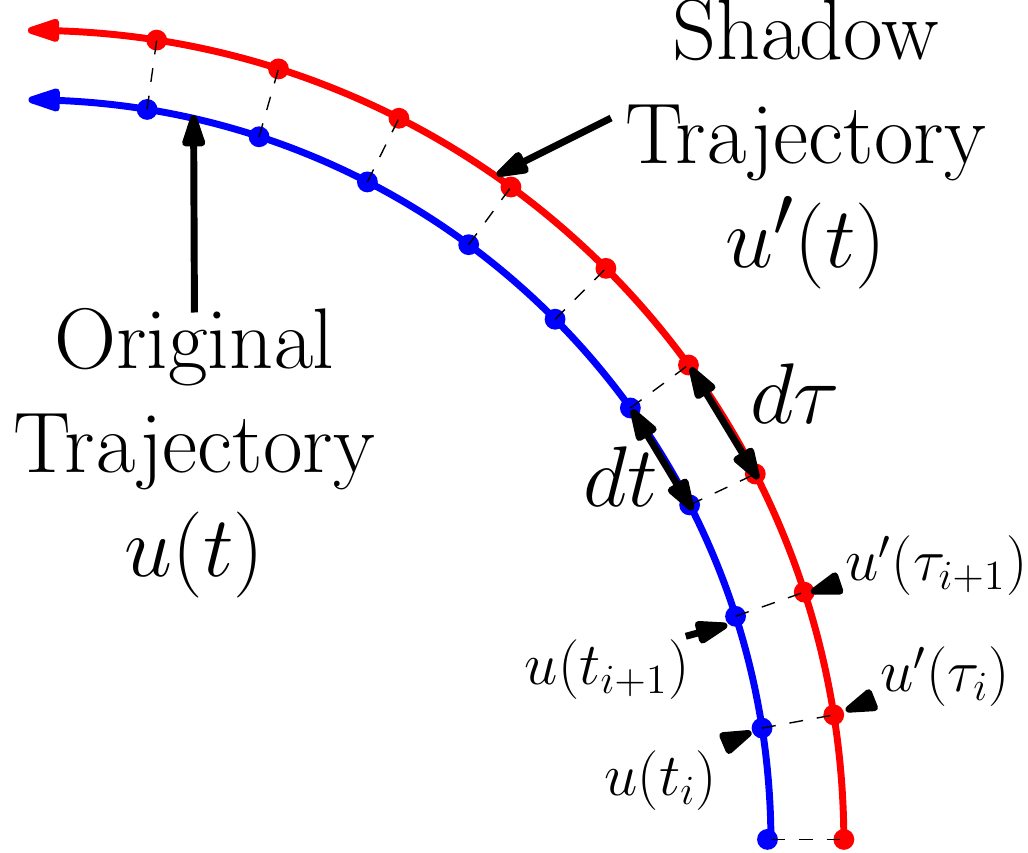}
\caption{LEFT: Original and shadow phase space trajectories without any time transformation ($d\tau/dt = 1$).  RIGHT: Original and shadow phase space trajectories with a time transformation $d\tau/dt = 1 + \eta$ that minimizes the distance between the two trajectories in phase space for all time.  Note that $\eta$ is a ``time stretching factor'' and that $-1.0<\eta<1.0$.  }
\label{f:time_dilation}
\end{figure}

The time transformation alluded to in the shadowing lemma is required to deal with the zero (neutrally stable) Lyapunov exponent, whose covariant vector is simply $f(u)$.  This time transformation, referred to as ``time dilation'' in this paper and other LSS literature, is required to keep a phase space trajectory and its shadow trajectory close (in phase space) for all time as demonstrated in figure \ref{f:time_dilation}.  

\subsection{Computing the Shadow Trajectory}  

Although LSS could be implemented by computing Lyapunov exponents and covariant vectors (as in \cite{Wang:2013:LSS1}), this is not necessary.  To find the shadow trajectory, an optimization problem is solved, where the 
objective function is the L2 norm of the tangent solution. 
That is, for some system of equations $\frac{du}{dt} = f(u,\xi)$, 
the tangent equations are solved, where $v = \frac{\partial u}{\partial \xi}$:

\begin{equation}
 \min_{v,\eta} \frac{1}{2}\int_0^T v^2 + \alpha^2 \eta^2 dt, \quad s.t. \quad \frac{dv}{dt} = \frac{\partial f}{\partial u} v + \frac{\partial f}{\partial \xi} + \eta f \,, \quad 0<t<T ,
\label{e:opt_problem}
\end{equation}

\noindent where $\eta$ is the time dilation term, corresponding to the time transformation from the shadowing lemma discussed in the previous section.  For a detailed derivation of equation \eqref{e:opt_problem}, see the appendices of ``Sensitivity computation of chaotic limit cycle oscillations.'' by Q. Wang, R. Hui and P. Blonigan, available at arXiv:1204.0159.  

This optimization problem is a linearly constrained 
least-squares problem, with the following KKT equations, 
derived using calculus of variations (see appendix \ref{s:KKTderive}):

\begin{align}
\frac{\partial w}{\partial t} &= -\left(\frac{\partial f}{\partial u}\right)^* w - v \quad
w(0) = w(T) = 0 \label{e:KKTw}\\
\alpha^2 \eta &= - \langle f, w \rangle \label{e:KKTeta}\\
\frac{dv}{dt} &= \frac{\partial f}{\partial u} v + \frac{\partial f}{\partial s} + \eta f \label{e:KKTv}
\end{align}

Equations (\ref{e:KKTw}), (\ref{e:KKTeta}) and (\ref{e:KKTv}) can be combined to form a single second order equation for the Lagrange multiplier $w$, as in \cite{Blonigan:2013:LSS}:

\begin{gather}
-\frac{d^2 w}{d t^2} - \left(\frac{d}{d t} \left(\frac{\partial f}{\partial u}\right)^*  - \frac{\partial f}{\partial u} \frac{d}{d t} \right) w + \left(\frac{\partial f}{\partial u}\left(\frac{\partial f}{\partial u}\right)^* +\frac{1}{\alpha^2} ff^*\right) w = \frac{\partial f}{\partial \xi} \nonumber\\ w(0) = w(T) = 0\label{e:SPD_LSS}
\end{gather}

Equation (\ref{e:SPD_LSS}) shows that the LSS method has changed the tangent
equation from an initial value problem in time to a boundary value problem in time. The LSS solution can be used to compute gradients for some quantity of interest $J$ (i.e. Drag):  
\begin{equation}
\frac{\partial \bar{J}}{\partial \xi} = \overline{\left\langle \frac{\partial J}{\partial u}, v \right\rangle} + \overline{\eta J} - \overline{\eta} \overline{J}
\label{e:LSSgrad}
\end{equation}
\noindent Where $\overline{x} \equiv \frac{1}{T} \int_0^T x \ dt$. See appendix \ref{s:LSSgrad} for a derivation of equation (\ref{e:LSSgrad}).   

\subsection{Solving the KKT system numerically}
\label{s:discKKT}

Equations (\ref{e:KKTw}), (\ref{e:KKTeta}) and (\ref{e:KKTv}) are discretized using finite differences and combined to form the following symmetric system:  

{\scriptsize
\[
\left( \begin{array}{ccccc|cccc|cccc}
\rowcolor{red!70} I & & & & & & & & & F_0^T & & & \\
\rowcolor{red!70} & I & & & & & & & & G_1^T & F_1^T & & \\
\rowcolor{red!70} & & \ddots & & & & &  & & & G_2^T & \ddots & \\
\rowcolor{red!70} & & & I & & & & & &   & & \ddots & F_{m-1}^T \\
\rowcolor{red!70} & & & & I & & & & & & & & G_m^T \\\hline
 & & & & & \alpha^2 & & & & f_1^T &  & &  \\
 & & & & & & \alpha^2 & & &  & f_2^T &  &  \\
 & & & & & & & \ddots & & & & \ddots & \\
 & & & & & & & & \alpha^2 &  & & & f_m^T \\\hline
\rowcolor{yellow!70} F_0 & G_1 & & & & f_1 & & & & & & & \\
\rowcolor{yellow!70} & F_1 & G_2 & & & & f_2 & & & & & & \\
\rowcolor{yellow!70} & & \ddots & \ddots & & & & \ddots & & & & & \\
\rowcolor{yellow!70} & & & F_{m-1} & G_m & & & &  f_m & & & & \end{array} \right) \left(\begin{array}{c} 
\rowcolor{red!30} v_0\\
\rowcolor{red!30} v_1 \\
\rowcolor{red!30} \vdots \\
\rowcolor{red!30} \vdots \\
\rowcolor{red!30} v_m \\\hline
 \eta_1 \\
 \eta_2 \\
 \vdots \\
 \eta_m \\\hline
\rowcolor{yellow!30} w_1 \\
\rowcolor{yellow!30} w_2 \\
\rowcolor{yellow!30} \vdots \\
\rowcolor{yellow!30}  w_m \end{array} \right) = -\left(\begin{array}{c} 
\rowcolor{red!70} 0 \\
\rowcolor{red!70} 0 \\
\rowcolor{red!70} \vdots \\
\rowcolor{red!70} \vdots \\
\rowcolor{red!70} 0 \\\hline
 0 \\
  \\
  \vdots\\
 0 \\\hline
\rowcolor{yellow!70} b_1\\
\rowcolor{yellow!70} b_2 \\
\rowcolor{yellow!70} \vdots \\
\rowcolor{yellow!70} b_m \end{array} \right)
\]}
\begin{gather*}
 F_i = \frac{I}{\Delta t} + \frac{1}{2}\frac{\partial f}{\partial u}(u_i,\xi), \quad G_i = -\frac{I}{\Delta t} + \frac{1}{2}\frac{\partial f}{\partial u}(u_i,\xi) ,\\ \quad b_i = \frac{1}{2} \left( \frac{\partial f}{\partial \xi}(u_i,\xi) + \frac{\partial f}{\partial \xi}(u_{i+1},\xi)\right), \quad f_i = \frac{1}{2} \left( f(u_i) + f(u_{i+1}) \right), \quad i = 0,...,m
\end{gather*}

This KKT system is a block matrix system, where each block ($I$, $G_i$, and $F_i$) is $n$ by $n$, 
where $n$ is the number of states (i.e. the product of the number of nodes and flow variables (density, velocities, enthalpy) in a Computational Fluid Dynamics (CFD) 
simulation).  $w_i$ and $v_i$ are length $n$ vectors, and $\eta_i$ is a scalar.  The blocks highlighted in red correspond to equation (\ref{e:KKTw}), white to equation (\ref{e:KKTeta}), and yellow to (\ref{e:KKTv}).  

Note that as the system is symmetric, the adjoint is computed simply 
by changing the right hand side, allowing many gradients to be computed 
simultaneously \cite{Giles:2000:adj}.  

\begin{figure}[!ht]
\begin{minipage}[b]{0.475\linewidth}
 \centering
\includegraphics[width=\textwidth]{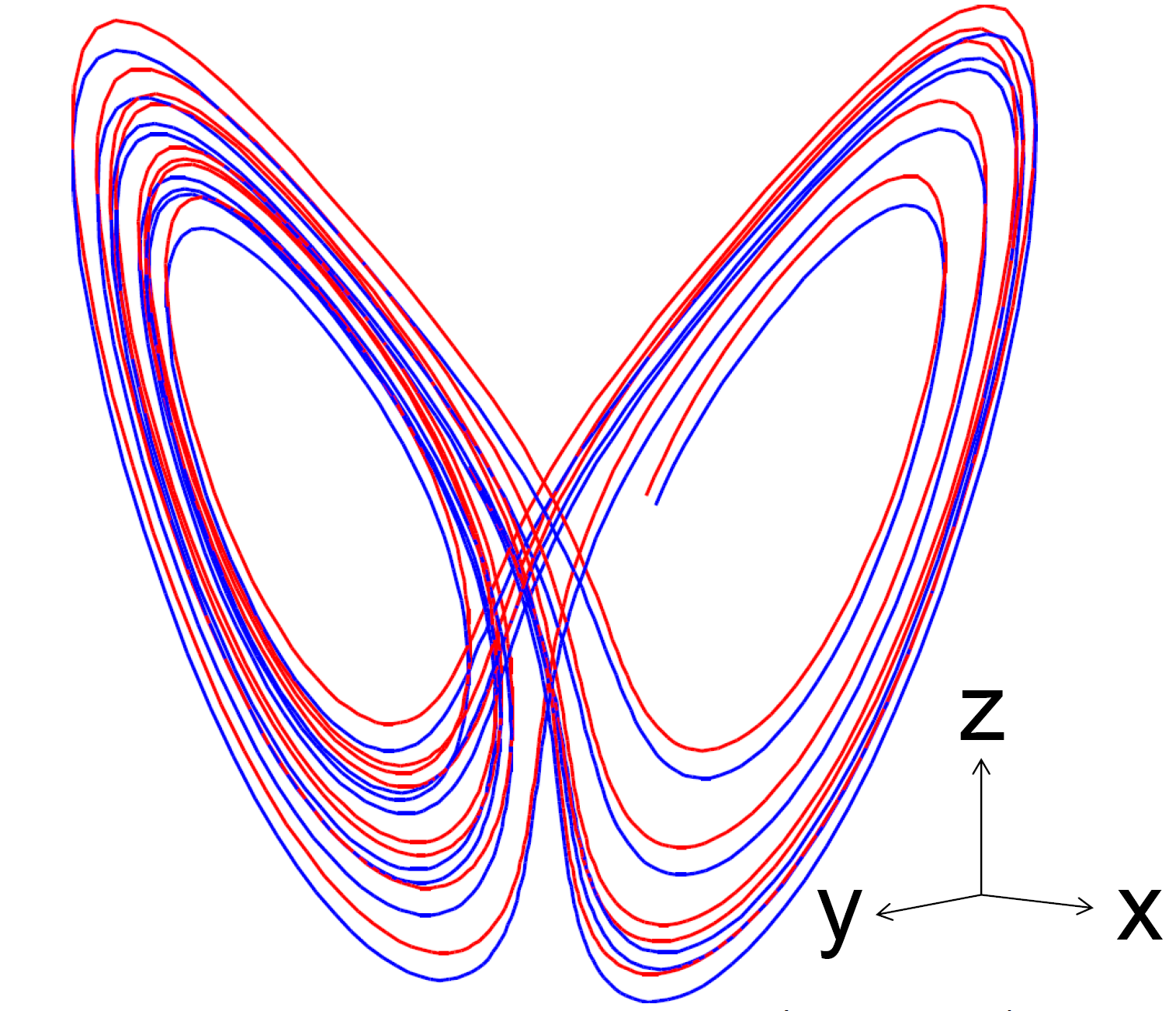}
\caption{Lorenz equation phase space trajectory ($u(t)$) for $r = 28$ (blue) and a corresponding approximate shadow trajectory ($u(t)+v(t)$) (red).  Integration time was $T=20$ in dimensionless time units.  }
\label{f:lorenz_shadow}
\end{minipage}
\hspace{0.125in}
\begin{minipage}[b]{0.475\linewidth}
 \centering
 \includegraphics[width=\textwidth]{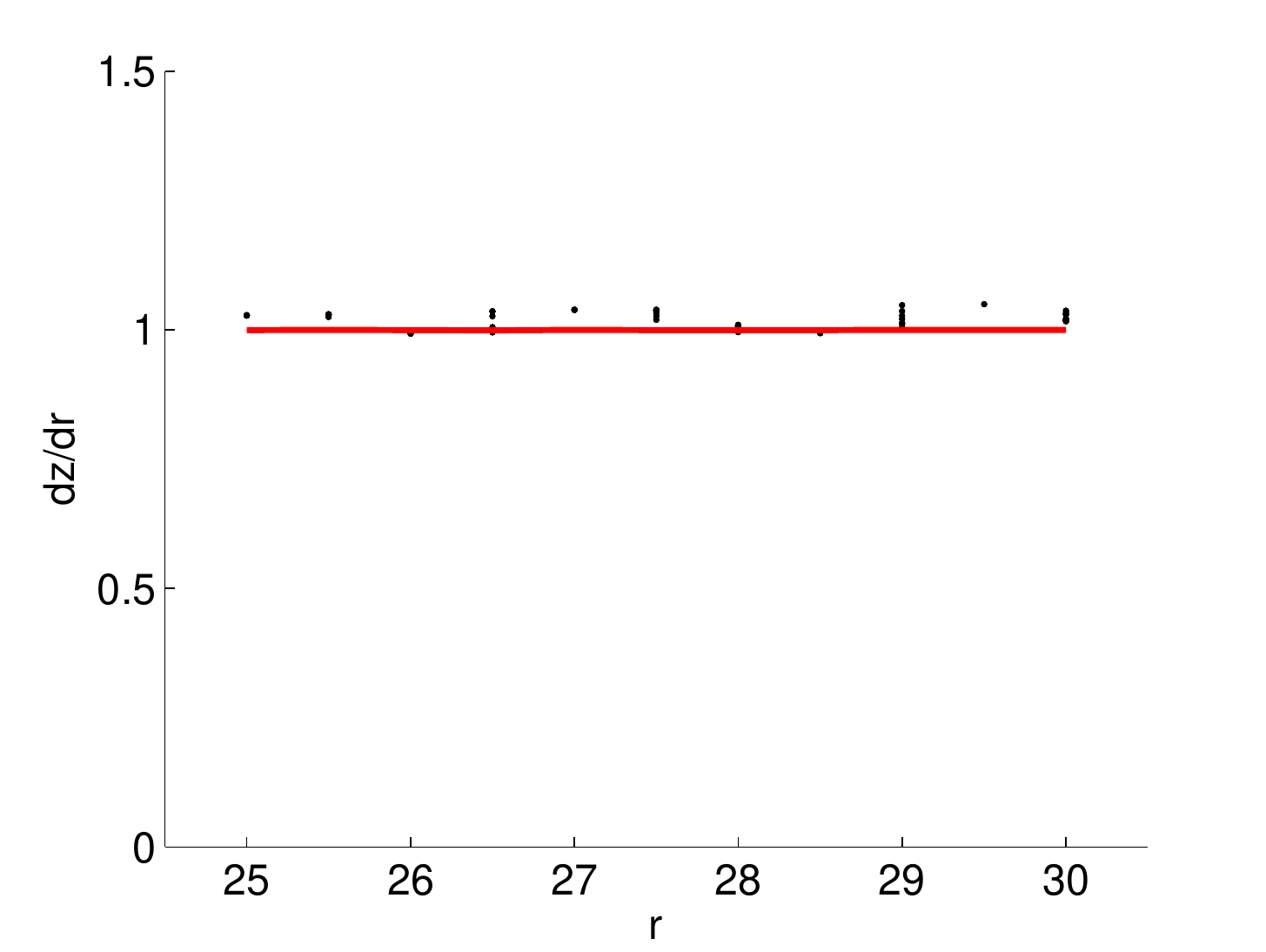}
 \caption{Gradient of long-time averaged $z$ with respect to the parameter $r$.  Gradients computed with trajectory length $T=20$ are shown as black diamonds.  Those computed using $T=1000$ are shown as a red line.  }
 \label{f:lorenz_grads}
\end{minipage}
\end{figure}

For an example of using LSS, we consider the Lorenz system.  The Lorenz system is a low order model of Rayleigh-B\'{e}nard convection \cite{Lorenz:1963:det}:

 \begin{equation}
 \frac{dx}{dt} = s (y - x), \quad
 \frac{dy}{dt} = x (r - z) - y, \quad
 \frac{dz}{dt} = x\,y - b\,z\;.
 \label{e:lorenz}
 \end{equation}

The Lorenz system was solved forward in time using a 4th order Runge-Kutta time stepping scheme.  To solve the KKT system for the Lorenz system, the KKT system is rearranged to a banded block matrix system, as shown in ``Sensitivity computation of chaotic limit cycle oscillations.'' by Q. Wang, R. Hui and P. Blonigan.  This system is then solved using a direct method.  Solutions were integrated from a random initial condition and run for 100 time units before LSS was applied, to ensure that the portion of the solution being used was on the attractor manifold.  

Figure \ref{f:lorenz_shadow} shows an approximate shadow trajectory for the Lorenz system.  Note how the trajectory stays close for all times, as governed by the shadowing lemma.  The shadow trajectory allows the computation of accurate sensitivities for relatively short integration times as shown in figure \ref{f:lorenz_grads}. The gradients shown in figure \ref{f:lorenz_grads} are within the bounds of the gradient found using linear regression: $1.01\pm0.04^*$ \cite{Wang:2013:LSS1}.  

Although rearranging the system to a banded system and solving with a direct works well for smaller systems like the Lorenz system, it would not be very efficient for LSS with larger systems.  

Firstly, the KKT system is quite large, with $2mn + n + 1$ by $2mn + n + 1$ 
elements for $m$ time steps.  For a discretization with a stencil of five 
elements, the matrix would have approximately $23mn$ non-zero elements.  
Consider a CFD simulation with $1\times10^5$ nodes and 
16000 time steps.  For a 2D compressible flow solver, there are four unknowns at each node, therefore, the number of degrees of freedom, $n$ is $4\times10^5$.  For this simulation, the KKT matrix would be $1.28\times10^{10}$ 
by $1.28\times10^{10}$ with $1.47\times10^{11}$ non-zero elements.  

Also, the bandwidth of the rearranged KKT system scales with the size of the system $n$, which means that direct methods, whose operation count scales with the matrix bandwidth squared times the rank, do not scale well for larger problems like CFD simulations.  On the other hand, an iterative method scales with the number of non-zero elements, not the product of bandwidth and rank.  Therefore, the operation count of an iterative method scales with $mn$, not $mn^3$ like a direct method.  
Because of this, iterative methods, such as multigrid, should be used to apply LSS to CFD simulations and other larger dynamical systems. 

\newcommand{\hw}{\underline{w}}
\newcommand{\hb}{\underline{b}}

\section{Multigrid-in-time for the Least Squares Shadowing Method}
\label{c:MG}

As discussed in section \ref{s:discKKT}, the KKT system is very large for many problems of interest, therefore we consider iterative methods to solve the KKT system.  As equation (\ref{e:SPD_LSS}) is a boundary value problem in time, a multigrid-in-time scheme is attractive because of its fast convergence relative to other iterative methods for many boundary value problems \cite{Briggs:2000:MGtutorial}. 

Several multigrid schemes were used to solve the LSS KKT system, with vastly varying convergence rates and computational costs. The first scheme considered, referred to as ``classic'' multigrid, is based upon the geometric multigrid scheme outlined by Briggs et al. \cite{Briggs:2000:MGtutorial}.  This scheme leads to very slow convergence of the KKT system's residual.  The second scheme considered was cyclic reduction.  This scheme converges in one cycle if the system is solved exactly on the coarsest grid.  However, the implementation of this scheme would use too much memory or require too many floating point operations to be viable for solving large scale systems with the computational resources that are currently available to most engineers and scientists.  Finally, it was found that using a scheme with a Krylov subspace solver as a smoother and higher order averaging between the KKT systems on the fine and coarse grids led to textbook multigrid convergence rates. This was found to be the case for schemes which satisfied the variational condition and those that did not, although schemes that did not satisfy the variational condition converged slightly slower than those that did.  

This section is comprised of a detailed discussion of the above results. Each subsection will outline each scheme.  Also, each subsection will present the performance of each scheme when used to apply the LSS method to the Lorenz equations.  Finally, the implications of these results for the computational efficiency of each scheme will be discussed.  

\subsection{KKT system Schur complement}

Rather than directly discretizing equation \eqref{e:SPD_LSS}, we solve the Schur complement of the KKT system shown in section \ref{s:discKKT}.  The KKT system can be written as:

\[
\left( \begin{array}{ccc}
I & 0 & B^T \\
0 & \alpha^2 I & C^T \\
B & C & 0 
\end{array} \right)
\left( \begin{array}{c}
\underline{v} \\
\underline{\eta} \\
\underline{w}
\end{array} \right)
= -\left( \begin{array}{c}
0 \\
0 \\
\underline{b}
\end{array} \right)
\]

\noindent where $B$ is a $mn\times (m+1)n$ matrix and $C$ is a $mn\times m$ matrix.  Conducting block Gaussian elimination, the Schur complement is found to be:

\[
 \underbrace{(BB^T + \frac{1}{\alpha^2} CC^T)}_{A}\underline{w} = \underline{b}
\]

\noindent Written in terms of the block matrices in section \ref{s:discKKT}:
{\small
 \begin{equation}
 \underbrace{\left( \begin{array}{ccc}
F_0 F_0^T + G_1 G_1^T + \frac{1}{\alpha^2} f_1 f_1^T & G_1 F_1^T &  \\
F_1 G_1^T & F_1 F_1^T + G_2 G_2^T +\frac{1}{\alpha^2} f_2 f_2^T  & G_2 F_2^T  \\
\ddots & \ddots & \ddots  \\
 & F_{m-1} G_{m-1}^T & F_{m-1} F_{m-1}^T + G_m G_m^T + \frac{1}{\alpha^2} f_m f_m^T  
 \end{array} \right)}_{A}
 \label{e:schurKKT}
\end{equation}
}
From this form we see that the Schur complement is a $mn \times mn$ SPD and block tridiagonal matrix.  

\subsection{Classic Multigrid}
\label{ss:CMG}

The first multigrid scheme implemented was a simple geometric scheme, referred to in this paper as ``classic multigrid''.  Injection is used for restriction and prolongation is carried out by linear interpolation. Injection is also used to restrict the system of equations \eqref{e:schurKKT}.  A ``V'' cycle was used, in which the system is coarsened until only one equation remains, then prolongated back to the full fine grid.  Block Gauss-Seidel iterations (one block per time step) were used for relaxation, with 4-10 cycles before restriction and after prolongation on each level.  

The scheme is found to be unstable if the under-relaxation factor for the Gauss-Seidel iterations is kept the same on the coarser grid.  An empirical formula was used to reduce the under-relaxation factor as the time-step $\Delta t$ increases in length on the coarse grid.  

The method was tested on the Lorenz system.  The gradient of time-averaged $z$ with respect to the parameters $b$, $r$, and $s$ was computed:
 
 \begin{equation}
 \frac{d\bar{z}}{ds} = 0.122, \quad \frac{d\bar{z}}{dr} = 1.00, \quad \frac{d\bar{z}}{db} = -1.67 .
\label{e:lorenz}
\end{equation}

The gradients with respect to $r$ and $b$ are within the error bounds of the gradients obtained by \cite{Wang:2013:LSS1} using linear regression, while the gradient with respect to $s$ was slightly over-predicted.
 
\begin{figure}
\begin{minipage}[b]{0.475\linewidth}
 \centering
 \includegraphics[width=\textwidth]{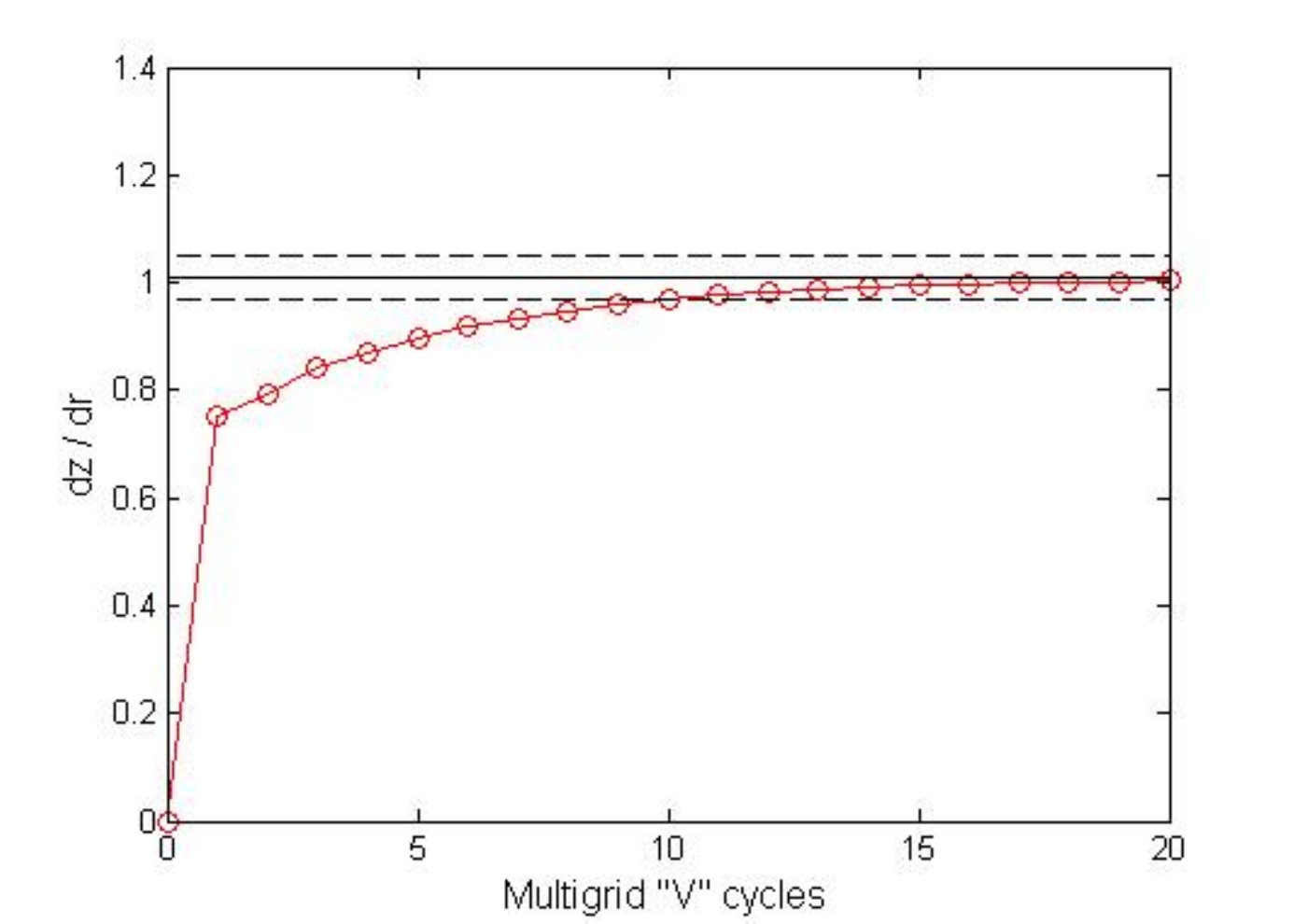}
 \caption{Convergence of the gradient of time-averaged z with respect to $r$, as computed using multigrid in time with 10 relaxation iterations before restriction and after prolongation on each level.  }
 \label{f:grad_conv}
\end{minipage}
\hspace{0.04\linewidth}
\begin{minipage}[b]{0.475\linewidth}
 \centering
 \includegraphics[width=\textwidth]{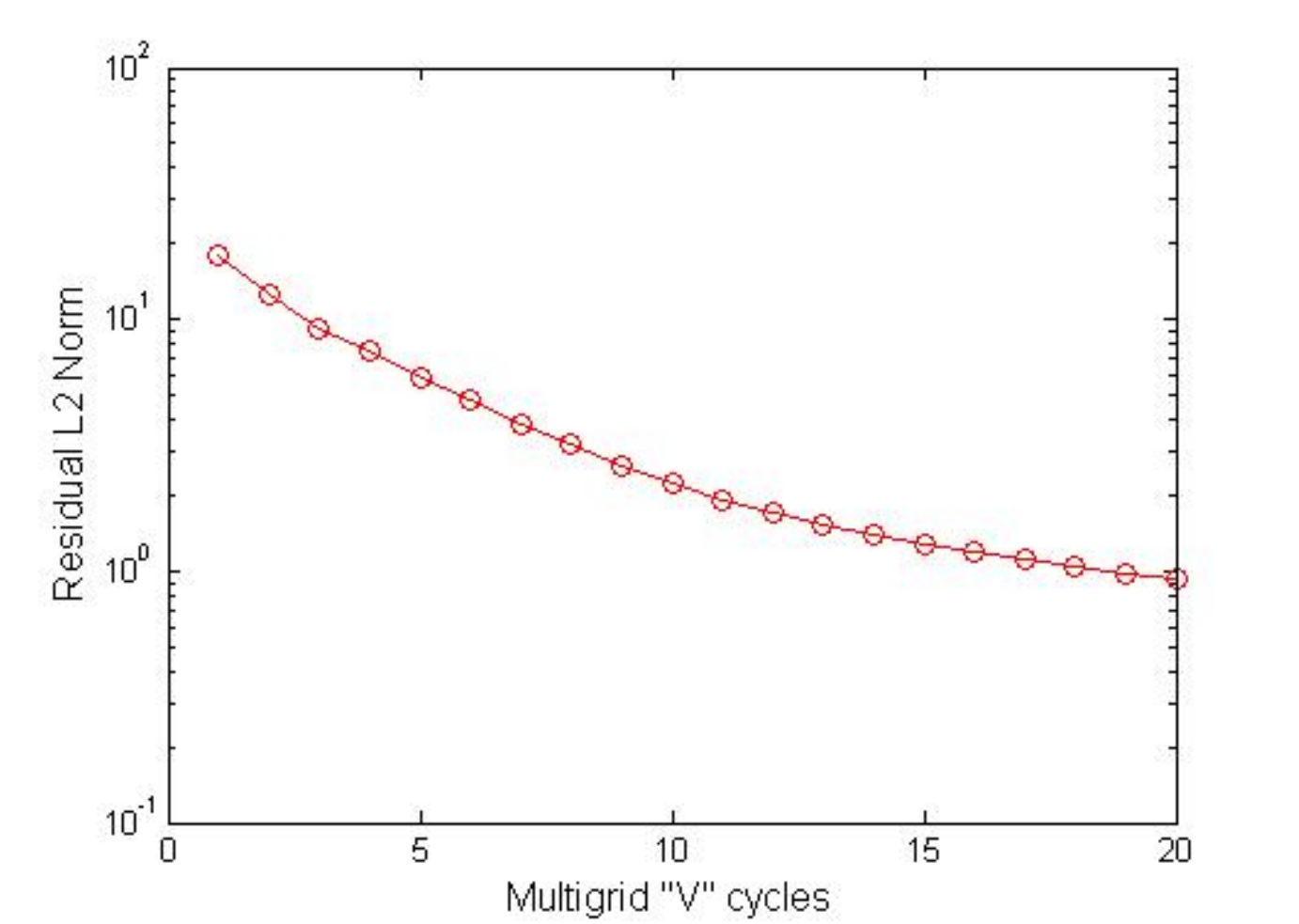}
 \caption{L2 norm of the residual while solving for the gradient of time-averaged z with respect to $r$ using multigrid in time.  Similar behavior was observed when computing other gradients}
 \label{f:res_conv}
\end{minipage}
\end{figure}

The gradients converged within 20-30 cycles, as shown in figure \ref{f:grad_conv}. However, the residual of the system did not converge as quickly, as in figure \ref{f:res_conv}.  In fact, the residual was not observed to converge to machine precision until around $10^4$ cycles.  Subsequent analysis explored the causes of the slow convergence of the residual. Firstly, the method was analyzed by conducting Ideal Coarse Grid (ICG) iterations, as defined in \cite{diskin:2004:MG}.  The convergence of the ICG iterations was found to be satisfactory, suggesting that the relaxation scheme was working well and the convergence issues arose from the grid coarsening scheme.  Also, convergence on individual grids was analyzed.  It was found that the residual did not decrease in magnitude on the coarsest grids.  
   
\noindent

Although this implementation of the multigrid in time method found correct gradient values, the slow convergence of the residual makes the robustness of the method questionable. The multigrid analysis methods carried out indicate that different grid coarsening techniques need to be explored to ensure textbook multigrid convergence to a near-zero residual.

\subsection{Cyclic Reduction}

In classic multigrid, linear interpolation is used for prolongation, but it is
not the best method when the coefficients of the equation being solved are highly oscillatory
or discontinuous \cite{Chan00}.  If this is the case, the coarse grid correction is a poor 
approximation for the low order error, and textbook multigrid convergence is not achieved.  
In this case, interpolation to some fine grid point with index $i$
should be carried out by solving (\ref{e:SPD_LSS}) with $\frac{\partial f}{\partial \xi}=0$ between
 the two nearest coarse grid points, $i-1$ and $i+1$, with the boundary conditions 
$w(t_{i-1}) = w_{i-1}$ and $w(t_{i+1}) = w_{i+1}$. It can be proved that this interpolation scheme 
leads to a multigrid scheme that converges independently of grid size for 1D problems \cite{Chan00}.  Furthermore, it can be shown that this scheme is equivalent to cyclic reduction \cite{Chan00}.  Cyclic reduction is carried out as follows: defining the lower, main and upper diagonal blocks in equation \eqref{e:schurKKT} as $L_i$, $D_i$ and $U_i$, elimination is conducted as follows for three given rows:

\begin{equation}
 L_{i-1} w_{i-2} + D_{i-1} w_{i-1} + U_{i-1} w_i = b_{i-1}
\label{e:row1}
\end{equation}
\begin{equation}
 L_i w_{i-1} + D_i w_i+ U_i w_ {i+1} = b_i
\label{e:row2}
\end{equation}
\begin{equation}
 L_{i+1} w_i + D_{i+1} w_{i+1} + U_{i+1} w_{i+2} = b_{i+1}
\label{e:row3}
\end{equation}

Equations (\ref{e:row1}) and (\ref{e:row3}) give expressions for $w_{i-1}$ and $w_{i+1}$, which can be substituted into equation (\ref{e:row2}):

\[
 L_I w_{i-2} + D_I w_i + U_I w_{i+2} = b_I
\]

\noindent Where:

\begin{equation}
\begin{array}{cc}
 L_I = & -L_i D_{i-1}^{-1} L_{i-1} \\
 D_I = & -L_i D_{i-1}^{-1} U_{i-1} + D_i - U_i D_{i+1}^{-1} L_{i+1} \\
 U_I = & -U_i D_{i+1}^{-1} U_{i+1} \\
 f_I = & -L_i D_{i-1}^{-1} f_{i-1} + f_i -U_i D_{i+1}^{-1} f_{i+1} \\
 b_I = & -L_i D_{i-1}^{-1} b_{i-1} + b_i -U_i D_{i+1}^{-1} b_{i+1}
\end{array}
\label{e:restriction}
\end{equation}

If the system is solved exactly on the coarsest grid, cyclic reduction will converge in one cycle.  Also, the algorithm can be implemented in parallel, as each coarse grid equation only depends on three adjacent fine grid equations.  

Although equation (\ref{e:restriction}) involves inverting the main diagonal matrices $D$, which contains products of Jacobians $\frac{\partial f}{\partial u}$, in practice these inversions do not need to be carried out (see appendix \ref{s:CRnoInvert}). However, the operation count of this scheme does not scale very well with the size of the KKT system.  It can be shown (see appendix \ref{s:CRflops}) that the number of floating point operations, $N$, required for matrix multiplication with the coarse grid matrix scales approximately as:

\[
N \sim \mathcal{O}(2p(2q)^l)
\]

Where $p$ is the number of operations required to multiply a vector by an $F_i$ or $G_i$ matrix, $q$ is the number of iterations needed by an iterative solver to solve any system involving $D_i^{-1}$, and $l$ is the number of levels, assuming the coarsest grid is $n \times n$ (one time step).  This is very large relative to the approximate operation count for a single Jacobi iteration on the LSS KKT system:

\[
N \sim \mathcal{O}(2^{(l+2)} p)
\]

To summarize, the advantages of cyclic reduction are in its memory efficiency and potential for parallel implementation, not its operation count.  

\subsection{Higher Order Averaging/Krylov Subspace Scheme}

Taking into account the performance of classic multigrid and cyclic reduction, a new multigrid scheme was designed for solving the KKT system associated with LSS.  The slow convergence of classic multigrid and the rapid convergence of cyclic reduction indicate that the accuracy of the coarse grid correction has a large effect on the convergence rate of the scheme.  To obtain a more accurate coarse grid correction, higher order averaging was used to coarsen the KKT system.  Higher order averaging ensures that the coarse grid non-linear solution $u(t)$ from which the KKT system is constructed is smooth, which leads to better multigrid performance \cite{Chan00}.  

Higher order averaging schemes are formed as follows: consider a two point (first order) average of some value $x$ at time step $i$:

\[
\overline{x}_i = \frac{1}{2} x_{i-1/2} + \frac{1}{2} x_{i+1/2}
\]

A first order scheme is formed by setting $x_{i-1/2} = x_{i-1}$ and $x_{i+1/2} = x_{i+1}$.  For a second order scheme, set $x_{i-1/2} = \frac{1}{2}x_{i-1} + \frac{1}{2}x_i$ and $x_{i+1/2} = \frac{1}{2}x_i + \frac{1}{2}x_{i+1}$:  

\[
\overline{x}_i = \frac{1}{4} x_{i-1} + \frac{1}{2} x_{i} + \frac{1}{4} x_{i+1}
\]

The third order scheme can be derived by substituting first order averages into the second order scheme, and so on:

\begin{align*}
\overline{x}_{i+1/2} = \frac{1}{8} x_{i-1} + \frac{3}{8} x_{i} + \frac{3}{8} x_{i+1} + \frac{1}{8} x_{i+2} \quad & \text{3nd Order} \\
\overline{x}_i = \frac{1}{16} x_{i-2} + \frac{1}{4} x_{i-1} + \frac{3}{8} x_{i} + \frac{1}{4} x_{i+1} + \frac{1}{16} x_{i+2} \quad & \text{4th Order} \\
\overline{x}_{i+1/2} = \frac{1}{32} x_{i-2} + \frac{5}{32} x_{i-1} + \frac{10}{32} x_{i} + \frac{10}{32} x_{i+1} + \frac{5}{32} x_{i+2} + \frac{1}{32} x_{i+3} \quad & \text{5th Order} 
\end{align*}

Higher order averaging was applied to multigrid-in-time in two ways.  {\em Matrix restriction multigrid} uses higher order averaging on the KKT matrix itself, using McCormick et al's variational conditions \cite{McCormick82}.  {\em Solution restriction multigrid} uses higher order averaging on the non-linear solution $u(t)$.  The KKT system is formed on the coarse grid from the restricted solution $u(t)$.  

In addition, it was observed that stationary iterative methods such as the Block Gauss Seidel method used in section \ref{ss:CMG} converge very slowly for the KKT system.  Krylov subspace methods such as conjugate gradient and MINRES\footnote{These methods were chosen because the KKT system is Symmetric Positive Definite.} were examined as smoothers as an alternative to the stationary methods used in classic multigrid.  Other parameters of the KKT system and the multigrid solver were also examined and some were found to have a considerable effect on convergence rates, especially the parameter $\alpha^2$ from equation (\ref{e:SPD_LSS}), the weighting of the time dilation term in the minimization statement.  

The following sections discuss matrix restriction multigrid and its application to LSS for the Lorenz equations, followed by a discussion of solution restriction multigrid.  The computational cost of both multigrid methods is compared to MINRES.  

\subsubsection{Matrix restriction multigrid}
\label{sss:matrixRes}

Matrix restriction multigrid was designed to satisfy the variational conditions 
\cite{McCormick82}:

\[
 I^{2h}_h = c^h (I^h_{2h})^T, \qquad A^{2h} = I^{2h}_h A^h I^h_{2h}
\]

Where $A^h$ is the fine grid matrix, $A^h$ is the coarse grid matrix, $I^{2h}_h$ is the
restriction matrix, $I^{h}_{2h}$ is the prolongation matrix and $c^h$ is some constant
that could depend on the grid dimension.  

Satisfying the variational condition ensures that the error of the solution will decrease 
monotonically, assuming the smoother decreases or does not change the magnitude of the error 
on all grids \cite{McCormick82}.

\begin{figure}
   \centering
   \begin{minipage}[b]{2.85in}
\includegraphics{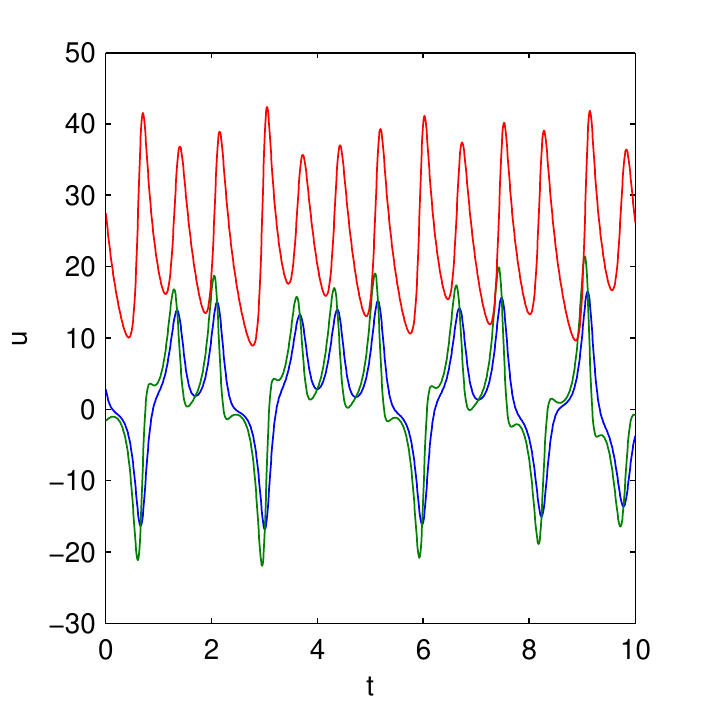}
\caption{Lorenz equation solution; x,y,z are blue, green and red respectively. \\}
\label{f:lorenz_soln}
   \end{minipage}
   \hspace{0.125 in}
   \begin{minipage}[b]{2.85in}
\includegraphics[width=2.8in]{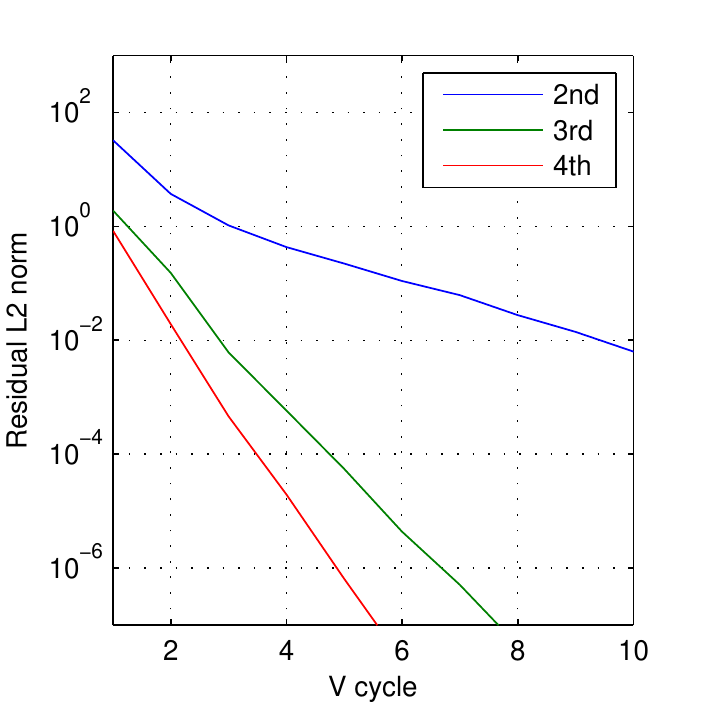}
\caption{Convergence of matrix restriction multigrid with $dt_f = 0.01$ for different orders of averaging.}
\label{f:interp_compare}   
   \end{minipage}
 \end{figure}

Matrix restriction multigrid has been demonstrated on LSS for the Lorenz equations, in particular the solution
shown in figure \ref{f:lorenz_soln}.  Unless otherwise stated, the results presented correspond to a multigrid
scheme with conjugate gradient (CG) smoothing with $\nu_1 = 30$ presmoothing iterations, 
$\nu_2 = 30$ postsmoothing iterations, 4th order averaging, $dt = dt_f = 0.0012$ on the
 finest grid and $dt = 0.08$ on the coarsest grid.  Although 30 is a large amount of smoothing iterations, it will be shown in section \ref{sss:SMG} that it results in a multigrid method that requires less operations than a Krylov subspace method on its own.  

Since the typical Lorenz equation solution in figure \ref{f:lorenz_soln} is oscillatory, a higher order averaging scheme can be used to improve the coarse grid correction by ensuring that the components of the block matrices (which depend on $u(t)$) that make up KKT matrix vary smoothly even on very coarse grids.  Because of this, 
the use of higher order averaging drastically improves the convergence rate of 
multigrid-in-time, as shown by figure \ref{f:interp_compare}. 

\begin{figure}
   \centering
   \begin{minipage}[b]{2.85in}
   \includegraphics{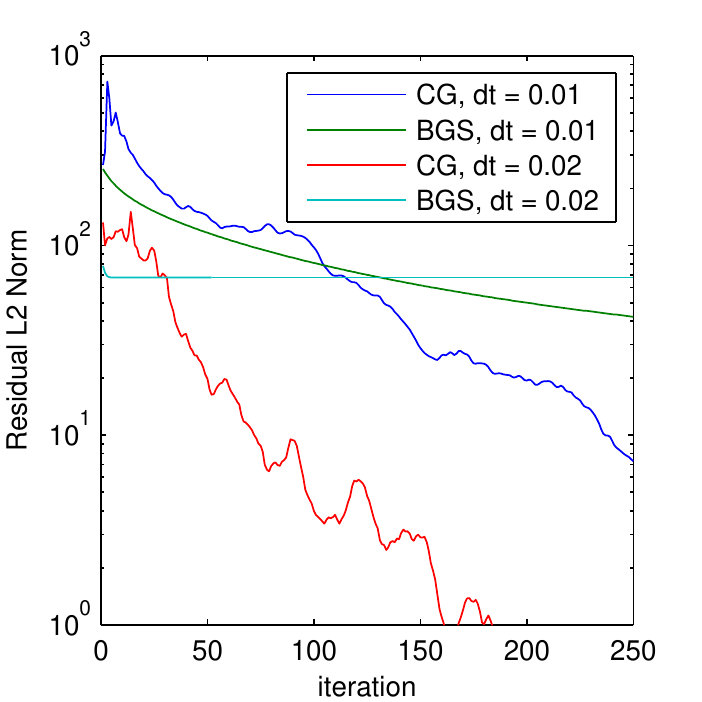}
\caption{Convergence of Block Gauss-Seidel and Conjugate Gradient on the fine grid 
(dt=0.01,$\alpha^2=40$) and a coarsened grid (dt=0.02,$\alpha^2=40$). \\ \\}
\label{f:GS_CG_coarse}
   \end{minipage}
   \hspace{0.125 in}
   \begin{minipage}[b]{2.85in}
\includegraphics[width=2.8in]{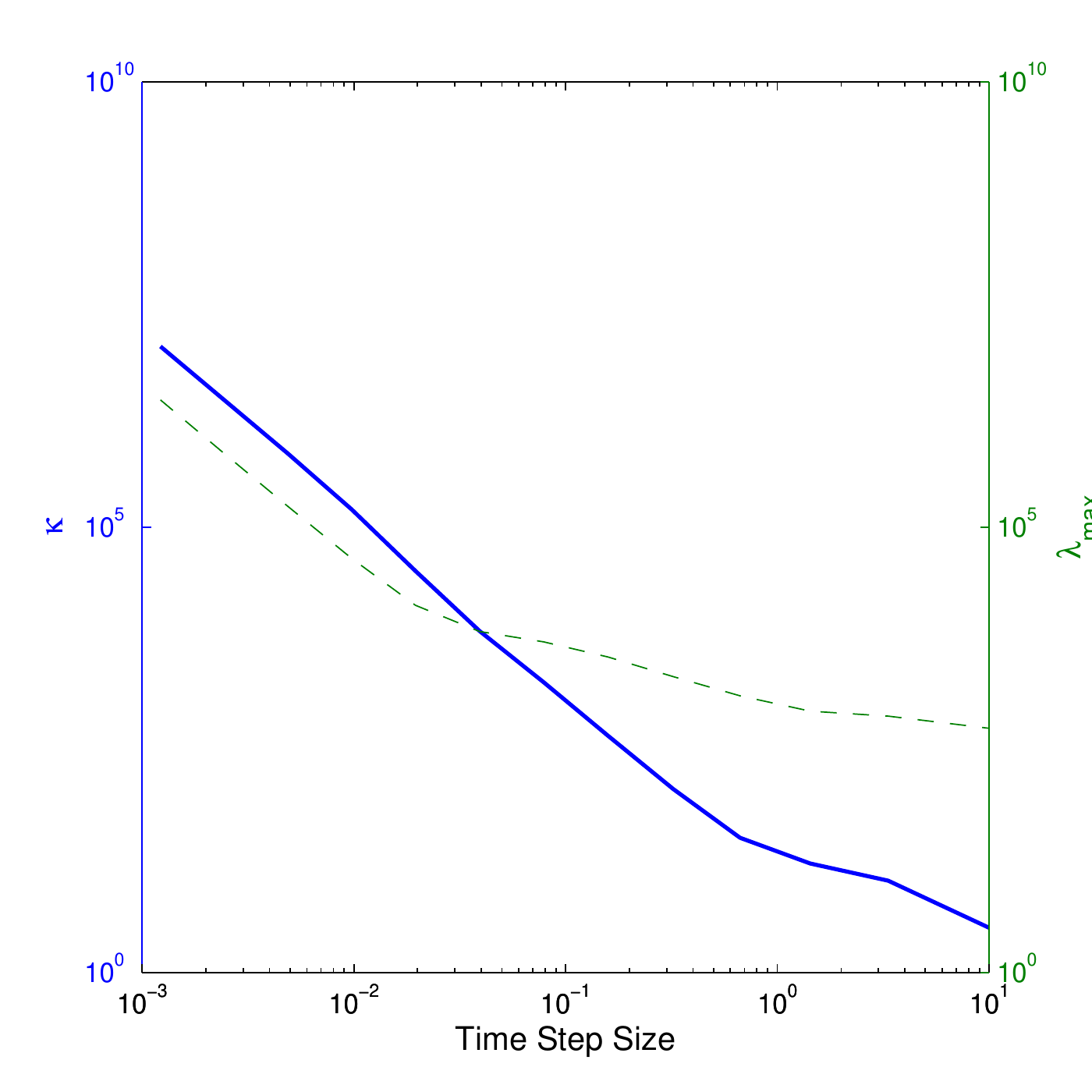}
\caption{Condition Number $\kappa$ (solid line) and maximum eigenvalue $\lambda_{max}$ (dashed line) versus time step $dt$ for coarsened grids corresponding
to a fine grid with $dt=0.001$ for the LSS system associated with the Lorenz equations}
\label{f:condition}
   \end{minipage}
 \end{figure}

Block Gauss-Seidel and other stationary solvers that are used for smoothing in classic multigrid schemes were found to converge very slowly.  Also, the convergence of these solvers worsen as the grid is coarsened, as
shown in figure \ref{f:GS_CG_coarse}.  This is because the largest eigenvalue of the KKT system, which can be shown to set the convergence rate of a stationary solver \cite{Golub1996}, decreases at a slower rate as the KKT system is coarsened, as shown in figure \ref{f:condition}.  This is in contrast to the behavior of the finite difference (or finite element) matrix of the Poisson equation, whose largest eigenvalue decreases exponentially as the grid is coarsened, resulting in much faster convergence of Jacobi or Gauss-Seidel solvers on the coarser grids \cite{Golub1996}.  

To accelerate convergence of multigrid-in-time, conjugate gradient (CG)
 is used as a smoother.  With a CG smoother, the rate of convergence is bounded by 
$((\kappa - 1)/\kappa)^N$, where $\kappa$ is the condition number of the KKT matrix $A_h$ and $N$ is the
number of smoothing iterations \cite{Douglas85}.  The condition number decreases quickly as the 
grid is coarsened (see figure \ref{f:condition}), leading to faster CG convergence on coarser
grids, as seen in figure \ref{f:GS_CG_coarse}.

\begin{figure}
   \centering

\includegraphics{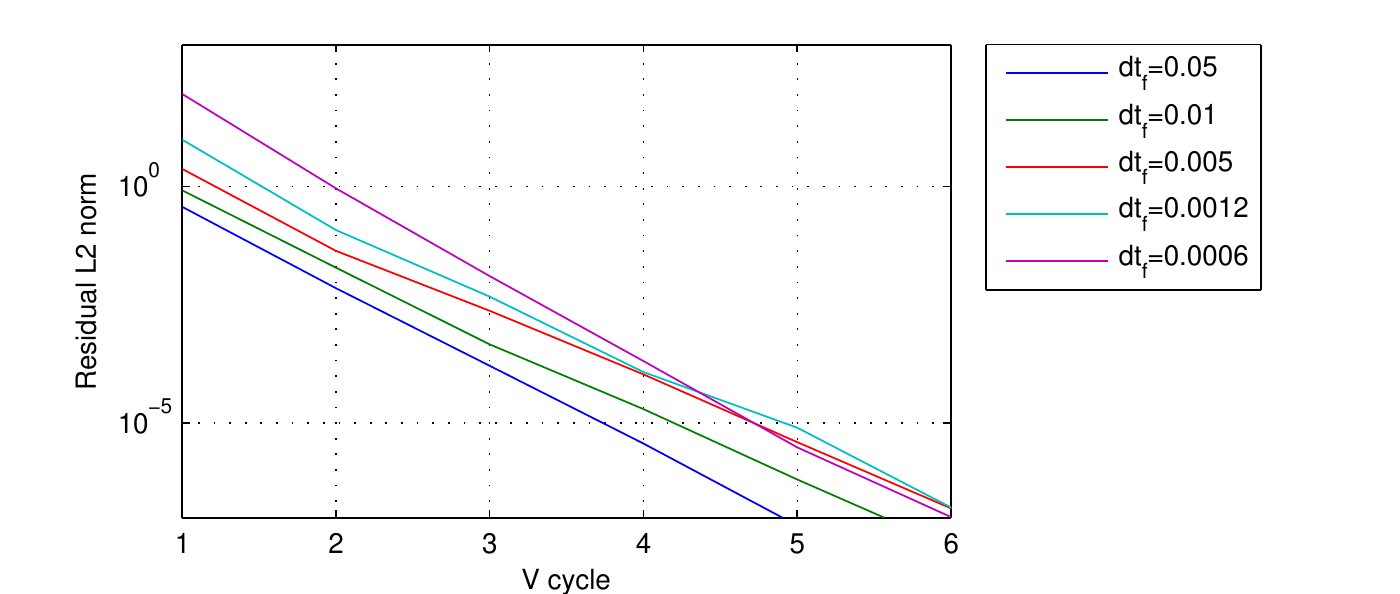}
\caption{Convergence of matrix restriction multigrid for different fine grid time steps $dt_f$. }
\label{f:dt_f_compare}

\includegraphics{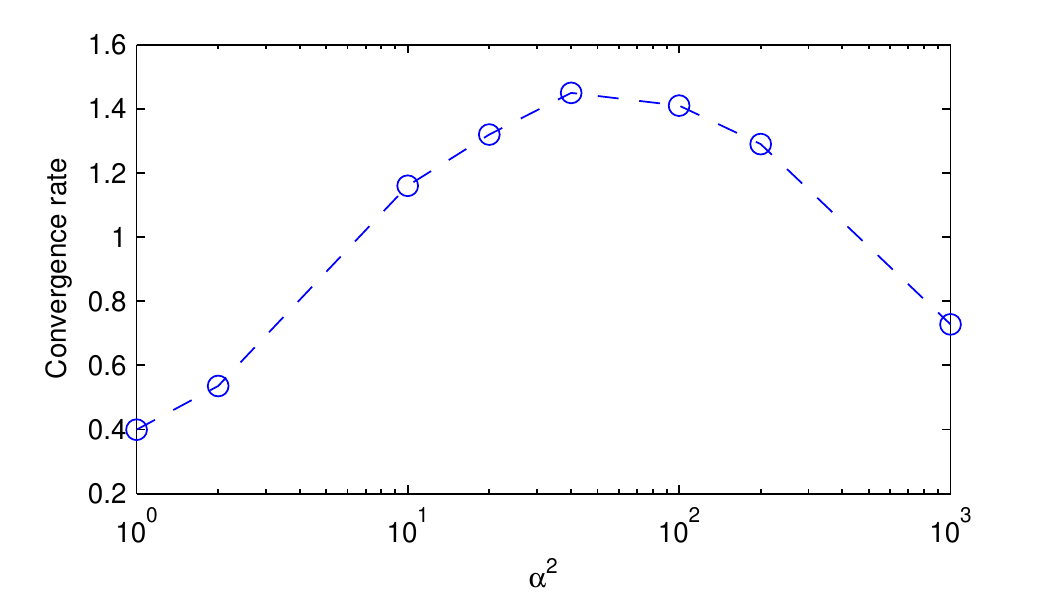}
\caption{Convergence rate of matrix restriction multigrid versus $\alpha^2$.  The convergence rate is $-\gamma$, from the curve fit $\log_{10} \|r\|_{L2} =  \gamma \log_{10} N_V + \log_{10} C$ of the residual L2 norm $\|r\|_{L2}$ versus the number of V-cycles $N_V$. } 
\label{f:alpha2_compare}

\includegraphics{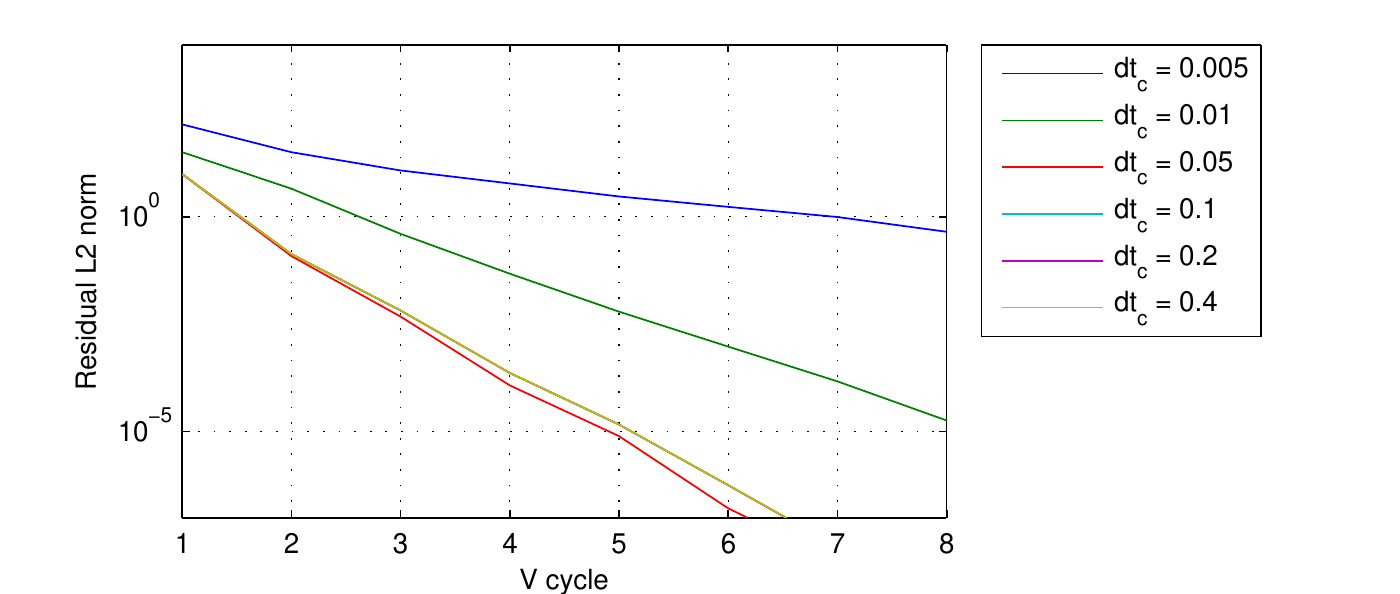}
\caption{Convergence of matrix restriction multigrid for different coarsening thresholds $dt_c$.  }
\label{f:dt_c_compare}

\end{figure}

A matrix restriction scheme with a CG smoother and 4th order averaging has been observed to
converge independently of the number of fine grid points, as shown in figure \ref{f:dt_f_compare}.

 In addition, a number of tuning parameters were observed to affect this convergence
rate. The parameter $\alpha^2$ from equation (\ref{e:SPD_LSS}), the weighting of the time dilation term in the
minimization statement, has a strong effect on multigrid convergence, as shown in figure
\ref{f:alpha2_compare}.  There is an optimal $\alpha^2$, found to be equal to around 40 in the
case of our sample problem, but this number is most likely specific to the Lorenz system. $\alpha^2$ has an effect on convergence because it affects the 
condition number and the eigenvalue spectrum of $A$ on all grids.  

The coarsening threshold $dt_c$ is defined as the time step size below which 
multigrid is not called recursively.  Figure \ref{f:dt_c_compare} shows that there is 
an optimal value for $dt_c$.  Below this, the coarse grid correction actually slows 
convergence in some cases, because it is a poor approximation for low order errors.   

\begin{figure}
	\centering
   \includegraphics{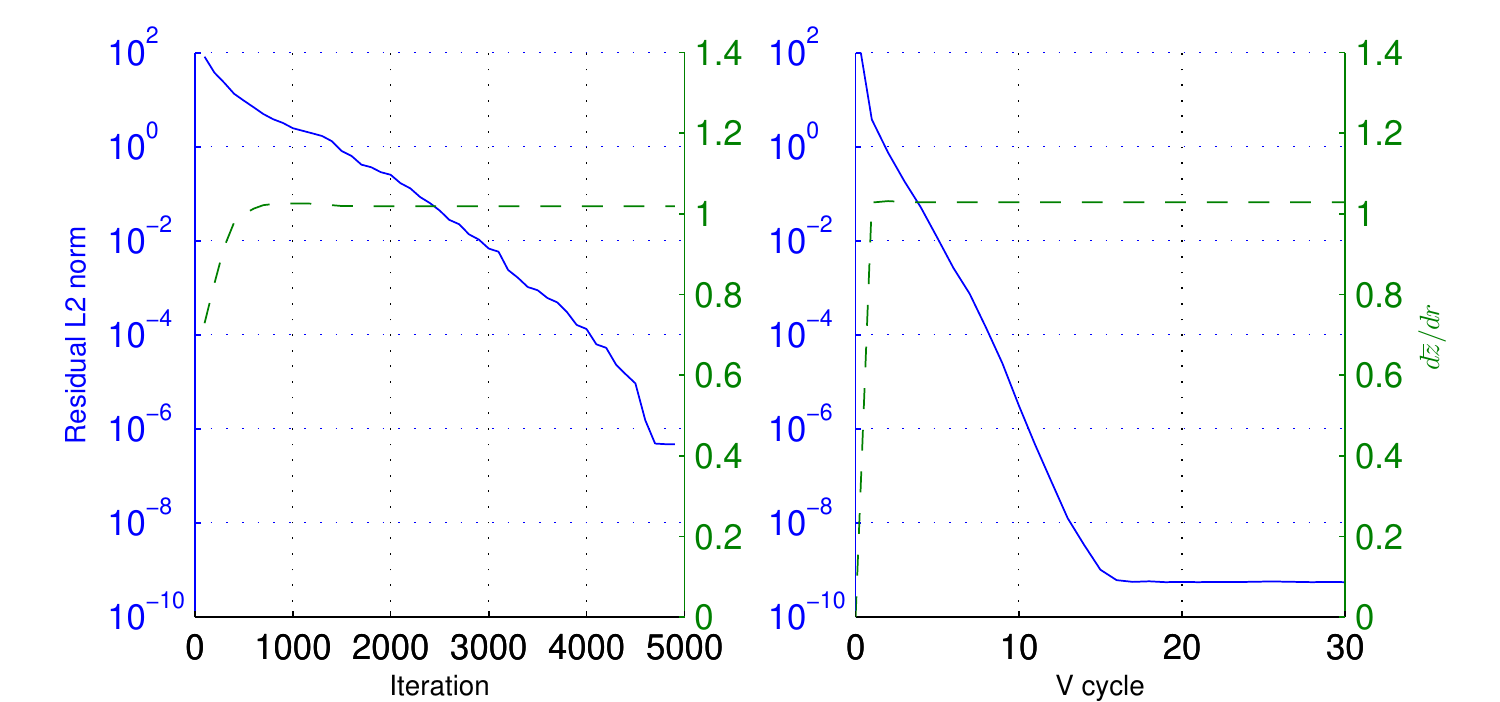}
\caption{LEFT: Convergence of MINRES for an LSS system for the Lorenz equations with $dt_f=0.004$ and $\alpha^2 = 40$.  RIGHT: Convergence of matrix restriction multigrid for an LSS system for the Lorenz equations.  The dashed lines shows the gradient computed at a given iteration, which should be roughly $1.01\pm0.04$ \cite{Wang:2013:LSS1} }
\label{f:grad_sol_mat}  
\end{figure}

To determine the relative efficiency of matrix restriction multigrid we compare its computational cost to that of using MINRES to solve the fine grid solution. Figure \ref{f:grad_sol} shows that the solution of the KKT system converges after about 4700 iterations, and the gradient converges after about 1800 iterations.  For a Krylov subspace method applied to a sparse matrix, the number of floating point operations for a single iteration, $p_{MINRES}$, is $p_{MINRES} \sim \mathcal{O}(mn)$,  $m$ is the number of time steps, $n$ is the number of dimensions of the dynamical system, and $mn$ is the size of the Schur complement of the KKT system in equation \eqref{e:schurKKT} \cite{Golub1996}.  Therefore, for the solution in figure \ref{f:grad_sol_mat} the total number of operations, $P_{MINRES}$, is roughly

\[
P_{MINRES} \sim 4700\mathcal{O}(mn)
\]

Using the variational condition to form the coarse grid matrix makes the matrix restriction multigrid more expensive than classic multigrid.  To conduct matrix multiplication on a coarse grid, the coarse grid solution is prolongated to the fine grid, multiplied by the fine grid matrix $A_h$ and then restricted to the coarse grid. Ignoring the cost of restriction and prolongation, for a fixed number of smoothing iterations $N = \nu_1 + \nu_2$ and 10 grids, the cost of one V-cycle of matrix restriction multigrid, $p_{MMG}$, is:

\begin{equation}
p_{MMG} \sim 10\mathcal{O}(mn(\nu_1 + \nu_2))) = 600 \mathcal{O}(mn)
\label{e:matrixMGcost}
\end{equation}

Figure \ref{f:grad_sol_mat} shows that the solution of the KKT system converges after about 17 cycles with matrix restriction multigrid, therefore:

\[
P_{MMG} \sim 10200 \mathcal{O}(mn) \approx 2P_{MINRES}
\]

Matrix restriction multigrid requires roughly twice as many operations as MINRES.  However, matrix restriction multigrid computes the correct gradient after 2 cycles, which requires $1200\mathcal{O}(mn)$ operations, slightly less than the roughly $1800\mathcal{O}(mn)$ operations required by MINRES.  

\subsubsection{Solution restriction multigrid}
\label{sss:SMG}

Although matrix restriction multigrid performs very well, it costs much more than MINRES.  We can reduce computational costs of multigrid by using solution restriction instead of matrix restriction.  Solution restriction only satisfies part of the variational condition, the restriction and prolongation operators are transposes of one another:  

\[
 I^{2h}_h = c^h (I^h_{2h})^T
\]

 \begin{figure}
   \centering
\includegraphics{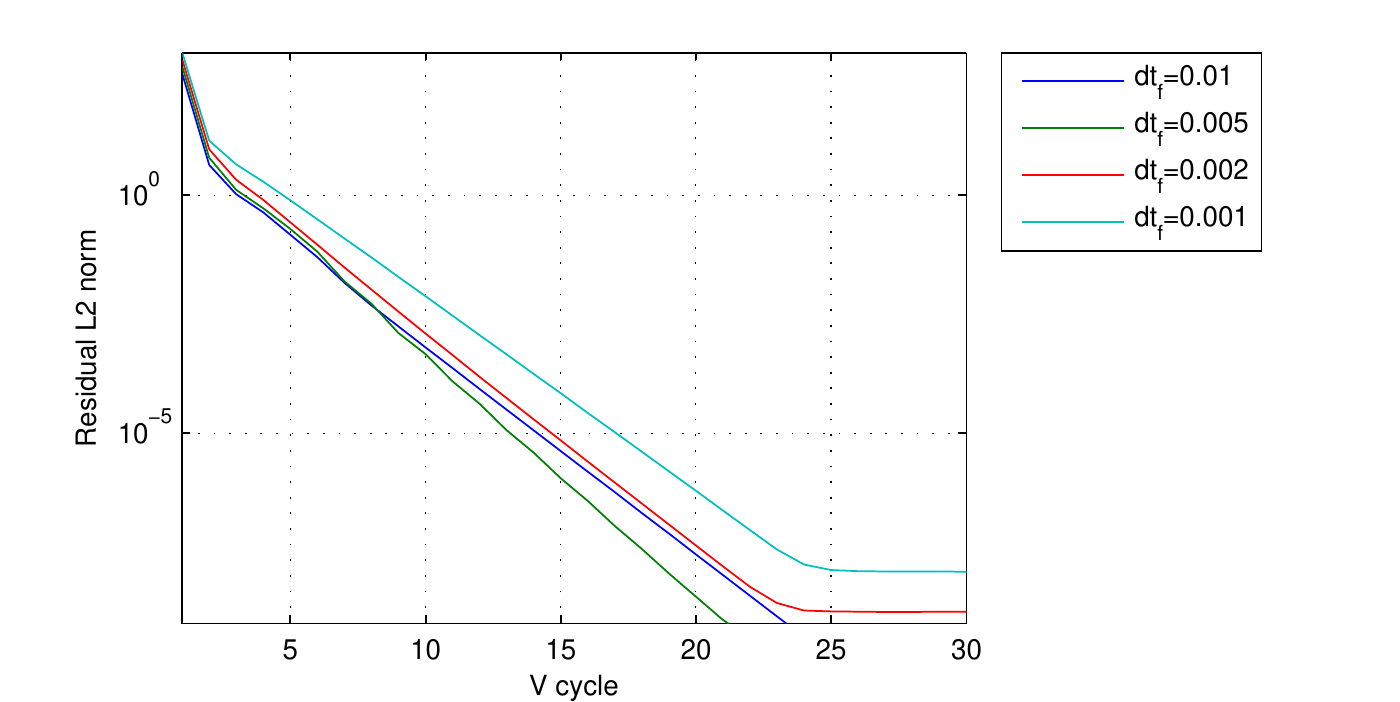}
\caption{Convergence of solution restriction multigrid for different values of $dt_f$.}
\label{f:sol_dtf_compare}
 \end{figure}

By restricting the non-linear solution $u(t)$ instead of the matrix $A_h$, solution restriction multigrid results in reduced smoothing costs on the coarse grid, because the KKT system formed on the coarse grid from the restricted $u(t)$ is half the size of the fine grid system. As shown by figure \ref{f:sol_dtf_compare}, the solution restriction scheme also leads to convergence rates with little dependence on $dt_f$ for a given $dt_c$ for the Lorenz equations \footnote{Parameter values used for this plots in this section (unless otherwise stated): $dt_f = 0.004$ or $m = 4096$,  $\nu_1=\nu_2 =30$, $dt_c = 0.2$, $\alpha^2 = 40$, MINRES smoothing, and 3rd order averaging}.

\begin{figure}
	\centering
   \includegraphics{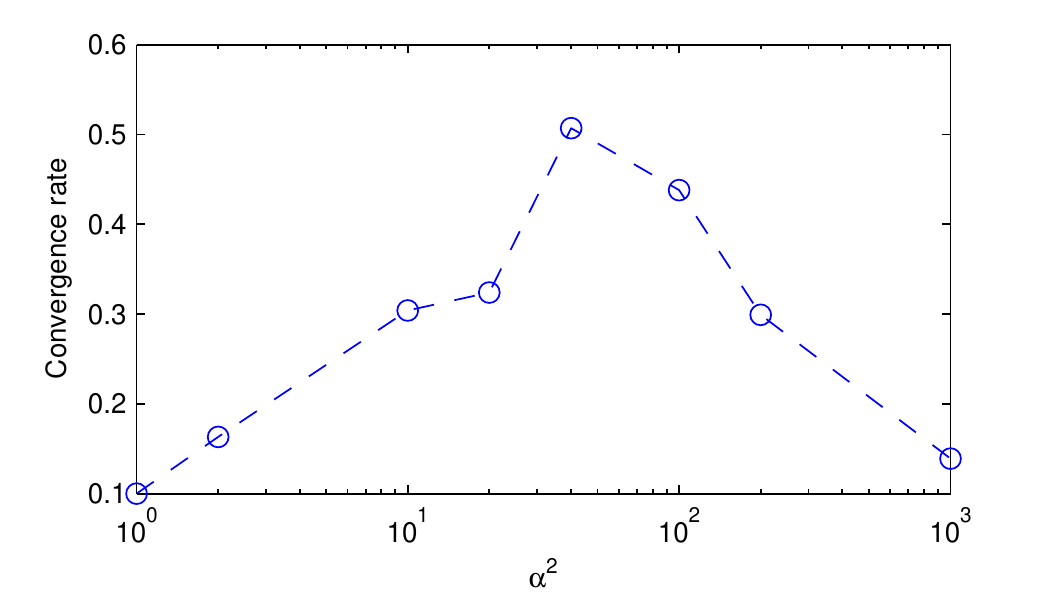}
\caption{Convergence rate of solution restriction multigrid versus $\alpha^2$.  The convergence rate is $-\gamma$, from the curve fit $\log_{10} \|r\|_{L2} =  \gamma \log_{10} N_V + \log_{10} C$ of the residual L2 norm $\|r\|_{L2}$ versus the number of V-cycles $N_V$. }
\label{f:sol_alpha2_compare}  
\end{figure}
   
\begin{figure}   
   \hspace{0.125 in}
   \begin{minipage}[b]{2.85in}
   \centering
   \includegraphics{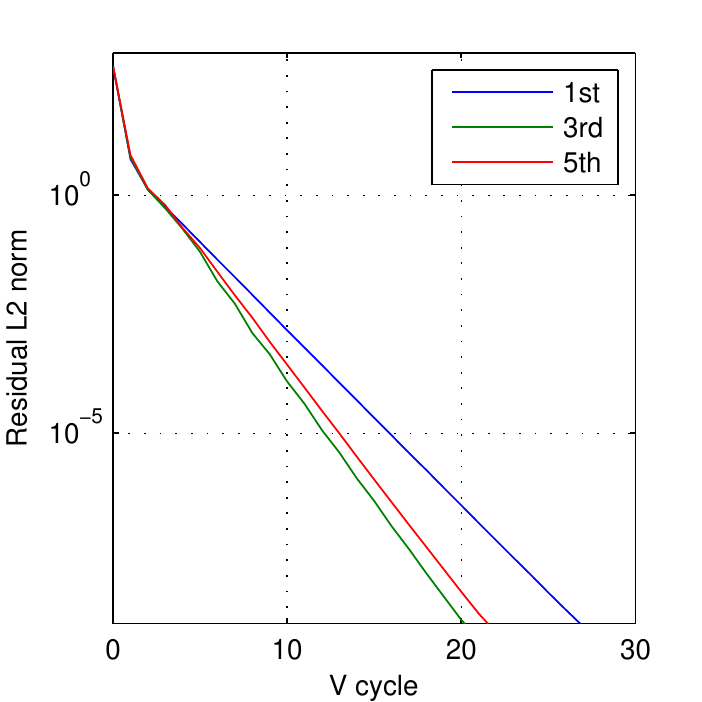}
\caption{Convergence of solution restriction multigrid with 1st, 3rd and 5th order averaging. \\ }
\label{f:sol_avg_compare}    
   \end{minipage}
   \hspace{0.125in}
   \begin{minipage}[b]{2.85in}
   \centering
   \includegraphics{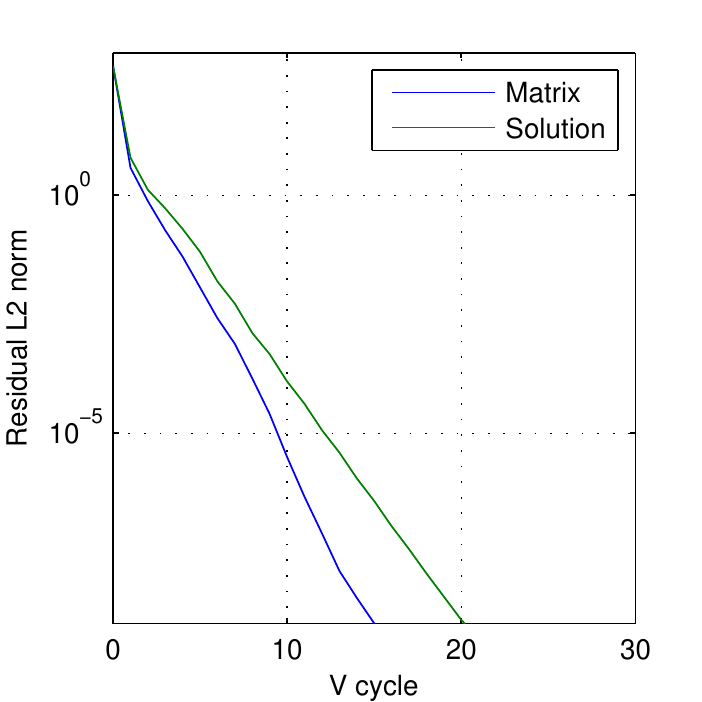}
\caption{Convergence plots for matrix restriction multigrid with MINRES smoothing, $dt_c = 5$, and 4th order averaging and solution restriction multigrid. }
\label{f:op_vs_sol}
   \end{minipage}
 \end{figure}

As observed for matrix restriction multigrid, the value of the parameter $\alpha^2$ and the order of averaging both affect the convergence rate, as shown in figures \ref{f:sol_alpha2_compare} and \ref{f:sol_avg_compare}.  However, higher order averaging is only beneficial to a certain degree, as the 5th order scheme leads to slower convergence than the 3rd order one (figure \ref{f:sol_avg_compare}).  This is because high order averaging could smooth $u(t)$ too much.  If this is the case the course grid solution is a poor approximation of the errors on the fine grid.  Also, $\alpha^2$ has a much greater effect than the averaging scheme on the convergence rate of solution restriction multigrid.  

There is a slight trade-off for the lower costs of solution restriction multigrid.  Figure \ref{f:op_vs_sol} shows that solution restriction multigrid converges slightly slower than matrix restriction multigrid. However, this slower convergence does not outweigh the benefits of the lower cost of solution restriction multigrid.

\begin{figure}
	\centering
   \includegraphics{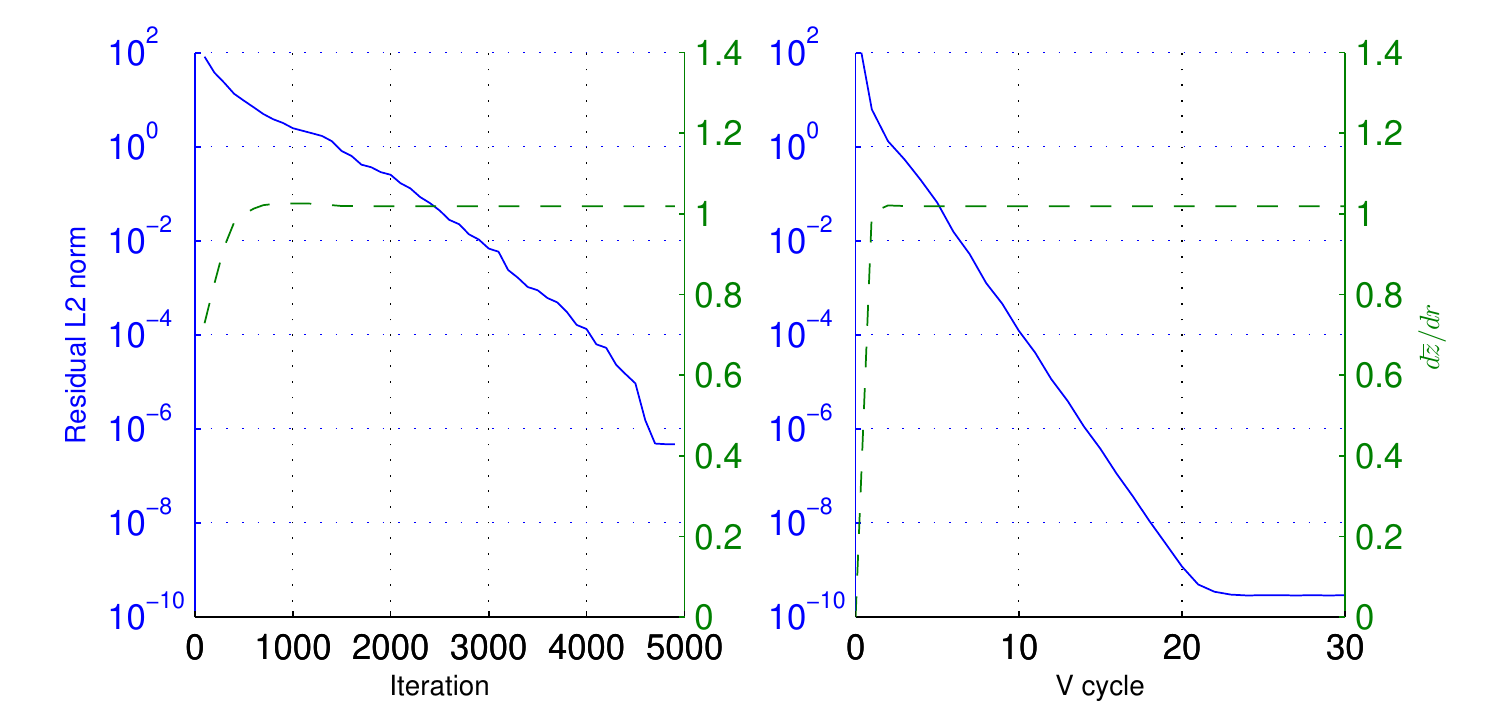}
\caption{LEFT: Convergence of MINRES for an LSS system for the Lorenz equations with $dt_f=0.004$ and $\alpha^2 = 40$.  RIGHT: Convergence of solution restriction multigrid for an LSS system for the Lorenz equations.  The dashed lines shows the gradient computed at a given iteration, which should be roughly $1.01\pm0.04$ \cite{Wang:2013:LSS1} }
\label{f:grad_sol}  
\end{figure}

These benefits can be seen by comparing the cost of solution restriction multigrid to that of matrix restriction multigrid and MINRES.  For a fixed number of iterations $N = \nu_1 + \nu_2$, the cost of smoothing is halved for one coarsening as $mn$ is halved when the system is coarsened.  When a large number of grids are used for a V-cycle, the cost of one V-cycle of solution restriction multigrid, $p_{SMG}$, is:

\begin{gather*}
p_{SMG} \sim \mathcal{O}(mn(\nu_1 + \nu_2)) + \frac{1}{2}\mathcal{O}(mn(\nu_1 + \nu_2)) + \frac{1}{4}\mathcal{O}(mn(\nu_1 + \nu_2)) + ... \\
 \approx 2 \mathcal{O}(mn(\nu_1 + \nu_2)) = 120\mathcal{O}(mn)
\end{gather*}

Figure \ref{f:grad_sol} shows that the solution of the KKT system converges after about 20 cycles with solution restriction multigrid, therefore:

\[
P_{SMG} \sim 2400 \mathcal{O}(mn) \approx \frac{1}{2} P_{MINRES} \approx \frac{1}{4} P_{MMG}
\]

In addition to requiring only half of the floating point operations of MINRES and a quarter of the operations of matrix restriction multigrid, solution restriction multigrid computes the correct gradient after 2 cycles, which requires $240\mathcal{O}(mn)$ operations, an order of magnitude less than the roughly $1800\mathcal{O}(mn)$ operations required by MINRES and a quarter of the cost of matrix restriction multigrid.

\section{Conclusion}
\label{c:conclusions}

In conclusion, a multigrid-in-time scheme could be used to implement the least squares shadowing (LSS) method in a relatively computationally efficient manner.  To solve the large KKT system associated with LSS, several multigrid-in-time schemes have been investigated.  Classic geometric multigrid with a Gauss-Seidel smoother was found to converge very slowly.  Cyclic reduction converged in one cycle and can be run in parallel but it used a very large number of floating point operations to solve the system if the number of dimensions of the system, $n$, was very large. A higher order averaging multigrid scheme with a Krylov subspace smoother was found to give textbook multigrid convergence.  Matrix restriction multigrid converges quickly, but smoothing on the coarse grid requires a similar amount of operations to smoothing on the fine grid.  Solution restriction multigrid converges slightly slower than matrix restriction multigrid, but is considerably less expensive.  It was also found that the parameter $\alpha^2$, the weighting of the time dilation term in equation (\ref{e:opt_problem}), had a large effect on the rate of convergence of matrix restriction and solution restriction multigrid.  

Because of its lower cost, solution restriction multigrid is currently the most promising numerical method for implementing LSS for large scale systems. Therefore, future work will start with further investigation of solution restriction multigrid-in-time.  In particular, a method to determine the $\alpha^2$ value corresponding to the fastest rate of convergence needs to be found.  

Once a robust, scalable numerical method for solving the KKT equation is found, LSS will be extensively tested on chaotic and turbulent fluid flows, such as homogeneous, isotropic turbulence and turbulent channel flow.  Eventually, LSS could be used to investigate more complicated flows such as flow around a lifting body, the fuel injection system in a jet engine or scramjet combustor, or an internal flow in a rocket engine. 

\section*{Acknowledgements}
The authors would like to acknowledge AFSOR Award F11B-T06-0007 under Dr. Fariba Fahroo, NASA Award NNH11ZEA001N under Dr. Harold Atkins as well as financial support from the NDSEG fellowship. 
 
\pagebreak

\appendix
\section{Least Squares Shadowing Sensitivity Analysis}

\subsection{Deriving the KKT System}
\label{s:KKTderive}

\noindent First form the Lagrangian function for equation (\ref{e:opt_problem}):
\[
L = \int_0^T \frac{v^2 + \alpha^2 \eta^2}{2} + \left\langle w,\left( -\frac{dv}{dt} + \frac{\partial f}{\partial u} v + \frac{\partial f}{\partial \xi} + \eta f \right) \right\rangle \ dt
\]

\noindent Now consider the first variation, which can be rearranged and transformed by integration by parts:
\[
\delta L = \int_0^T  \left\langle \delta v , \left( v + \frac{dw}{dt} + \frac{\partial f}{\partial u}^* w \right) \right\rangle \ dt + \int_0^T (\alpha^2 \eta + \langle f, w \rangle)\delta \eta \ dt + \left. \langle \delta v, w \rangle \right|_0^T 
\]

\noindent For first order optimality $\delta L = 0$ for all $\delta v$ and $\delta \eta$, therefore: 

\begin{equation}
\boxed{
v + \frac{dw}{dt} + \frac{\partial f}{\partial u}^* w = 0 , \quad w(0) = w(T) = 0
}
\label{e:lag_mult}
\end{equation}

\noindent and

\begin{equation}
\boxed{
\alpha^2 \eta + \langle f, w \rangle = 0
}
\label{e:time_dilation}
\end{equation}

\noindent Equations (\ref{e:lag_mult}) and (\ref{e:time_dilation}), along with the tangent equation, make up the KKT system.

\subsection{Computing Sensitivities using a Shadow Trajectory}
\label{s:LSSgrad}

\noindent For a trajectory $u(t)$ and shadow trajectory $u'(t)$:

\begin{align*}
\delta \bar{J} &= \frac{1}{\mathcal{T}} \int_0^{\mathcal{T}} J(u'(\tau)) \ d\tau - \frac{1}{T} \int_0^T J(u(t)) \ dt \\
&=  \int_0^{T} \frac{1}{\mathcal{T}} J(u'(\tau(t))) \frac{d\tau}{dt} + \frac{1}{T} J(u(t)) \ dt 
\end{align*}

\noindent For some perturbation to the tangent equation $\delta f = \varepsilon \frac{\partial f}{\partial \xi}$, the time transformation is $\frac{d\tau}{dt} = 1 + \varepsilon \eta $:

\[
\mathcal{T} = \int_0^T \frac{d\tau}{dt} dt = T + \varepsilon \underbrace{\int_0^T \eta \ dt}_{\equiv H} = T + \varepsilon H
\]

\noindent Therefore:
\[
\delta \bar{J} = \frac{1}{T+\varepsilon H} \int_0^{T} J(u'(\tau(t))) - J(u(t)) +   \varepsilon \eta J(u'(\tau(t))) - \frac{\varepsilon H}{T} J(u(t)) \ dt 
\]

\noindent Diving through by $\varepsilon$:

\begin{align*}
\frac{\partial \bar{J}}{\partial \xi} &= \lim_{\varepsilon \to 0} \left[\frac{1}{T+\varepsilon H} \int_0^{T} \frac{(J(u'(\tau(t))) - J(u(t))}{\varepsilon} +  \eta J(u'(\tau(t))) - \frac{H}{T} J(u(t)) \ dt \right] \\
&= \frac{1}{T} \int_0^T \left \langle \frac{\partial J}{\partial u}, v \right \rangle dt + \frac{1}{T} \int_0^T \eta J(u(t)) \ dt - \frac{H}{T} \frac{1}{T} \int_0^T J(u(t)) \ dt
\end{align*}

\noindent By the definition of $H$, and defining $\overline{x} \equiv \frac{1}{T} \int_0^T x \ dt$:

\[ \boxed{
\frac{\partial \bar{J}}{\partial \xi} = \overline{\left \langle \frac{\partial J}{\partial u}, v \right \rangle} + \overline{\eta J} - \overline{\eta} \overline{J} 
}\]

\section{Cyclic Reduction}
\subsection{Conducting cyclic reduction without inverting main diagonal matrices}
\label{s:CRnoInvert}

Consider the following system:

\[
 \left(\begin{array}{ccc}
D_1 & U_1 & 0 \\
L_2 & D_2 & U_2 \\
0 & L_3 & D_3 \end{array}\right) \left(\begin{array}{c}
w_1 \\
w_2 \\
w_3
\end{array}\right) = \left(\begin{array}{c}
 b_1 \\
b_2 \\
b_3
\end{array}\right)
\]

\noindent Applying equation (\ref{e:restriction}), the following system is obtained:

\[
 A w_2 = b
\]

\noindent with:

\[
 \begin{array}{cc}
 A = & -L_2 D_{1}^{-1} U_{1} + D_2 - U_2 D_{3}^{-1} L_{3} \\
 b = & -L_2 D_{1}^{-1} b_{1} + b_2 -U_2 D_{3}^{-1} b_{3}
\end{array}
\]

\noindent This system can be solved iteratively, using some preconditioner $P$:

\[
 P\Delta x = b - Ax_k,\quad x_{k+1} = x_k + \Delta x
\]

\noindent Where $x_k$ is the value of $w_2$ after $k$ iterations. To compute $A x_k$, it is decomposed into three parts:

\[
A x_k =  -L_2 D_{1}^{-1} U_{1} x_k + D_2 x_k - U_2 D_{3}^{-1} L_{3} x_k = \alpha + \beta + \gamma
\]

\noindent $\alpha$ and $\gamma$ include an inverted matrix, but the inversion can be avoided as follows:

\[
-L_2 D_{1}^{-1} U_{1} x_k = \alpha
\]

\noindent Compute $y_k = U_1 x_k$:
\[
-L_2 D_{1}^{-1} y_k = \alpha
\]

\noindent Next, define $z_k = D_{1}^{-1} y_k$.  Iteratively solve:

\[
 D_{1} z_k = y_k
\]

\noindent and use the result to compute $\alpha$:

\[
\alpha = -L_2 z_k
\]

$\gamma$ and the right hand side $b$ can be computed using a similar method.  This idea can be applied to a much larger system and allows cyclic reduction to be conducted without inverting any Jacobian matrices.

\subsection{Estimating the operation count for cyclic reduction}
\label{s:CRflops}

\noindent The following section estimates the operation count of the algorithm of section \ref{s:CRnoInvert}, in which block cyclic reduction is conducted without directly inverting any blocks of the matrix. First define: 
 
\begin{itemize}
\item $p$: the number of flops required to multiply a vector by a Jacobian matrix for the system of interest. 
\item $q$: the number of iterations required to carry out multiplication by an inverse matrix.  
\item $n$: the number of states in the system
\end{itemize}

As $L_i = F_{i-1} G_{i-1}^T$ and $U_i = G_i F_i^T $, multiplication by these matrices requires $2p$ flops. As $D_i = F_{i-1} F_{i-1}^T + G_i G_i^T + f_i f_i^T$, multiplication by these matrices requires $4p + 2n$ flops. For a Jacobi solver, the number of flops for inverse matrix multiplication is $qp$.  

An estimate of the operation count of cyclic reduction is shown for a few low order terms.  These were computed using a symbolic calculator.  Starting from a 1 by 1 coarse grid, the number of flops for multiplication by $D$, $U$ or $L$ were substituted into equation (\ref{e:restriction}) recursively.  The highest order term for a few different grids is shown in table \ref{t:CR_OP}, with the number of operations for a single Jacobi iteration (derived in the same way) for comparison. These estimates do not take into account fixed number of operations to backward substitute the coarse grid solution for the fine grid solution.  

\begin{table}[ht!]
\centering
\begin{tabular}[t]{|c|c|c|}\hline
Time steps $m$ & CR Flops & Jacobi Flops\\\hline
3 & $8pq^2$ & $20p + 13n$  \\
5 & $16pq^3$ & $36p + 23n$ \\
9 & $32pq^4$ & $68p + 43n$ \\
17 & $64pq^5$ &  $132p + 83n$ \\\hline
\end{tabular}
\caption{Estimate of Operation Count per iteration for cyclic reduction and for the Jacobi method for comparison.  Note that the cost of cyclic reduction is similar to that of a block Gaussian elimination scheme as block cyclic reduction is like applying block Gaussian elimination in parallel to a permuted system. }
\label{t:CR_OP}
\end{table}

\end{document}